\documentclass[12pt]{amsart}
\usepackage{amsmath,amssymb,bm}
\usepackage{graphicx,color}
\usepackage{subfigure}
\newtheorem{theorem}{Theorem}[section]
\newtheorem{lemma}[theorem]{Lemma}
\newtheorem{example}[theorem]{Example}
\newtheorem{remark}[theorem]{Remark}
\numberwithin{equation}{section}
\numberwithin{table}{section}
\numberwithin{figure}{section}
\allowdisplaybreaks[4]
\begin{document}
\def\O{\Omega}
\def\p{\partial}
\def\R{\mathbb{R}}
\def\bK{\mathbb{K}}
\def\argmin{\mathop{\rm argmin}}
\def\Hh{{H^{1/2}(\p\O)}}
\def\Har{\mathcal{H}(\O)}
\def\HOne{{H^1(\O)}}
\def\HOnez{{H^1_0(\O)}}
\def\LT{{L_2(\O)}}
\def\PH{P_{\mathcal{H}}}
\def\fC{\mathfrak{C}}
\def\LI{(-\Delta)^{-1}}
\def\LIh{(-\Delta_h)^{-1}}
\def\cT{\mathcal{T}}
\def\cA{\mathcal{A}}
\def\CPF{\mathrm{C}_{\mathrm{PF}}}
\def\fA{\mathfrak{A}}
\def\byh{\bar y_{h,\rho}}
\def\buh{\bar u_{h,\rho}}
\def\bph{\bar p_{h,\rho}}
\def\bps{\bar p_{*,\rho}}
\def\dyh{\dot y_h}
\def\dys{\dot y_{*}}
\def\duh{\dot u_h}
\def\duH{\dot u_{\sss H}}
\def\dph{\dot p_h}
\def\dps{\dot p_{*}}
\def\tyh{\tilde{y}_h}
\def\tuh{\tilde{u}_h}
\def\tys{\tilde{y}_*}
\def\tur{\tilde{u}_\rho}
\def\ss{\scriptstyle}
\def\sss{\scriptscriptstyle}
\def\MS{{V}_H^{\raise 2pt\hbox{$\ss \rm {ms},h$}}}
\def\MSLy{{y}_{\sss H,k}^{\mathrm{ms},h}}
\def\MSLv{{v}_{\sss H,k}^{\mathrm{ms},h}}
\def\MSLq{{q}_{\sss H,k}^{\mathrm{ms},h}}
\def\dMSLy{{\dot y}_{\sss H,k}^{\mathrm{ms},h}}
\def\dMSLp{{\dot p}_{\sss H,k}^{\mathrm{ms},h}}
\def\MSLz{{z}_{\sss H,k}^{\mathrm{ms},h}}
\def\bMSLy{{\bar y}_{\sss H,k}^{\mathrm{ms},h}}
\def\bMSLry{{\bar y}_{\sss H,k,\rho}^{\mathrm{ms},h}}
\def\bMSLp{{\bar p}_{\sss H,k}^{\mathrm{ms},h}}
\def\bMSLrp{{\bar p}_{\sss H,k,\rho}^{\mathrm{ms},h}}
\def\MSbK{\mathbb{K}_{\sss H,k}^{\mathrm{ms},h}}
\def\uH{u_{\sss H}}
\def\buH{{\bar u}_{\sss H}}
\def\MSV{V_{H,k}^{\hspace{1pt}\lower 3pt\hbox{$\ss\rm {ms},h$}}}
\def\PiH{\Pi_{\sss H}}
\def\KER{K_h^{\raise 1pt\hbox{$\ss \PiH$}}}
\def\cC{\mathcal{C}^{\raise 1pt\hbox{$\ss \PiH$}}}
\def\Ch{\cC_h}
\def\Chk{\cC_{h,k}}
\def\MS{{V}_H^{\raise 2pt\hbox{$\ss\rm {ms},h$}}}
\def\MSu{{u}_H^{\mathrm{ms},h}}
\def\cCk{\mathfrak{C}_{h,k}}
\def\red{\color{red}}
\def\bKsr{\bK_{*,\rho}}
\def\bys{\bar y_{*,\rho}}
\def\bur{\bar u_{*,\rho}}
%
%%%%%%%%%%%%%%%%%%%%%%%%%%%%%%%%%%%%%%%%%%%%%%%%%%%%%%
\title[A Multiscale FEM for an Optimal Control Problem with
 Rough Coefficients]
{A Multiscale Finite Element Method for an Elliptic Distributed
Optimal Control Problem \\ with Rough Coefficients and
Control Constraints}
 \author{Susanne C. Brenner}
\address{Department of Mathematics and Center for Computation
\& Technology,
Louisiana State University, Baton Rouge, LA 70803, USA}
\email{brenner@math.lsu.edu}
\author{Jos\'e C. Garay}
\address{Institute of Mathematics, Universt\"at Augsburg,
86159 Augsburg, Germany}
\email{jose.garay.fernandez@uni-a.de}
\author{Li-yeng Sung}
\address{Department of Mathematics and Center for Computation
\& Technology,
Louisiana State University, Baton Rouge, LA 70803, USA}
\email{sung@math.lsu.edu}
\begin{abstract}
  We construct and analyze a multiscale finite element method for an
  elliptic distributed  optimal control problem  with pointwise control constraints,
  where the state equation has rough coefficients.
  We show that the performance of the multiscale finite element method is similar to
  the performance of standard finite element methods for smooth problems and present
  corroborating numerical results.
\end{abstract}
\thanks{This work was supported in part
 by the National Science Foundation under Grant No.
 DMS-19-13035 and Grant No. DMS-22-08404.}
\date{September 11, 2023}
\maketitle
%%%%%%%%%%%%%%%%%%%%%%%%%%%%%%%%%%%%%%%%%%%%%%%%%%%%%%%%%
%%%%%%%%%%%%%%%%%%%%%%%%%%%%%%%%%%%%%%
\section{Introduction}\label{sec:Introduction}
 Let $\O\subset \R^d$  ($d=1,2,3$) be a polytopal domain, $y_d\in\LT$ and
 $\gamma\leq 1$ be a positive constant.
 The model optimal control problem (cf. \cite{Lions:1971:OC,Troltzsch:2010:OC})
  is to find
\begin{equation}\label{eq:OCP}
  (\bar y,\bar u)=\argmin_{(y,u)\in \bK} J(y,u),
\end{equation}
 where the cost function $J:\HOnez\times\LT\longrightarrow [0,\infty)$ is defined by
\begin{equation}\label{eq:JDef}
  J(y,u)=\frac12\big(\|y-y_d\|_\LT^2+\gamma \|u\|_\LT^2\big),
\end{equation}
 the closed convex subset $\bK$ of $\HOnez\times\LT$ is defined by the conditions
\begin{alignat}{3}
   a(y,z)&=\int_\O uz\,dx &\qquad&\forall\,z\in \HOnez,\label{eq:PDEConstraint}\\
        \phi_1\leq &\,  u\leq\phi_2&\qquad&\text{a.e. in $\O$},
        \label{eq:ControlConstraints}
\end{alignat}
 and the bilinear form  $a(\cdot,\cdot)$ is given by
\begin{equation}\label{eq:aDef}
  a(y,z)=\int_\O \cA\nabla y\cdot\nabla z\,dx.
\end{equation}
\par
 We assume that the components of the symmetric positive definite
 matrix $\cA$ belong to
 $L_\infty(\O)$, and that there exist positive constants $\alpha$ and $\beta$ such that
\begin{equation}\label{eq:Ellipticity}
\text{the eigenvalues of $\cA$ are bounded below (resp., above) by
 $\alpha$ (resp., $\beta$).}
\end{equation}
\par
  For the constraint functions $\phi_1$ and $\phi_2$, we assume
\begin{equation}\label{eq:ConstraintRegularity}
   \text{$\phi_1$ and $\phi_2$ belong to $H^1(\O)$},
\end{equation}
 and
\begin{equation}\label{eq:ConstraintInequality}
  \text{$\phi_1\leq\phi_2\;$ a.e. on $\O$.}
\end{equation}
\begin{remark}\label{rem:Notation}\rm
  Throughout this paper we follow the standard notation for differential operators,
  functions spaces and
  norms that can be found for example in
  \cite{Ciarlet:1978:FEM,ADAMS:2003:Sobolev,BScott:2008:FEM}.
\end{remark}
\begin{remark}\label{rem:RC}\rm
  The rough coefficients in the title of the paper refer to the fact that
  \eqref{eq:Ellipticity} is the only assumption on the matrix $\cA$.
  Under this assumption we have
  the relation
\begin{equation}\label{eq:SpaceComparison}
  \alpha|v|_{H^1(\O)}^2\leq \|v\|_a^2={a(v,v)}\leq \beta|v|_{H^1(\O)}^2
  \qquad\forall\,v\in H^1(\O)
\end{equation}
 and nothing more.  In particular, we do not assume the solution $y$ of
 \eqref{eq:PDEConstraint} belongs to
 $H^{1+s}(\O)$ for some positive $s$.
\end{remark}
\par
 It is well-known that standard finite element methods for elliptic boundary
 value problems with rough coefficients can
 converge arbitrarily slowly (cf. \cite{BO:2000:Bad}).  This is of course also
 the case for the optimal control
 problem defined by \eqref{eq:OCP}--\eqref{eq:aDef}.  Our goal is to design a
  multiscale finite element method
 whose performance is in some sense similar to that of the standard finite element
 methods for smooth problems.
\par
 The literature on the numerical solution of this optimal control problem
 is relatively small.
 For problems with scale separations and periodic structures, the method
 in \cite{LCY:2013:MS} is based on an asymptotic expansion of the solution,
  the method in
 \cite{CHLY:2015:MS} is based on the multiscale finite element space
  in \cite{CH:2003:Oscillating},
  and the method in \cite{GYWLY:2018:HMMOCP} is based on
 the heterogeneous multiscale method in \cite{EE:2003:HMSFEM}.
 For problems that do not
 assume scale separations or periodic structures,
 a numerical method based on the multiscale finite element space in
 \cite{CEL:2018:Multiscale}
 was investigated in \cite{AC:2022:DD26}, and a numerical method based on a
 generalization of the multiscale finite
 element space in \cite{OZB:2014:Homogenization}
 has just appeared in \cite{CLZZ:2023:MultiscaleCC}.
\par
 Our method is based on
 the local orthogonal decomposition (LOD) methodology
 (cf. \cite{MP:2021:LOD}) which,
 like the methods in \cite{AC:2022:DD26,CLZZ:2023:MultiscaleCC},
 does not require scale separations or periodic structures in the
 coefficient matrix $\cA(x)$.
 A variant of the LOD method for elliptic optimal control problems
 with rough coefficients but without
 control constraints can also be found in \cite{BGS:2022:LOD_OCP}.
\par
 The rest of the paper is organized as follows.  The properties of
 the continuous problem are recalled
 in Section~\ref{sec:Continuous} and  a  discretization of
 the optimal control problem is analyzed in
 Section~\ref{sec:FEMs}, where we present error estimates that are
 convenient for the
 error analysis of multiscale finite element methods.
 The construction and analysis of our multiscale finite element
 method are presented in Section~\ref{sec:DDLOD}, followed by
 numerical results in Section~\ref{sec:Numerics}.  We end with
 some concluding remarks
 in Section~\ref{sec:Conclusions}.
%%%%%%%%%%%%%%%%%%%%%%%%%%%%%%%%%%%%%%%%%%%%%%%%%%%%%%%%%%%%%%%%
\section{The Continuous Problem}\label{sec:Continuous}
 In this section we recall some well-known facts about the optimal
 control problem
 that can be found for example in
 \cite{Lions:1971:OC,Troltzsch:2010:OC}.
\par
  Since $\bK$ is nonempty under \eqref{eq:ConstraintInequality} and $J$
  is strictly
  convex and coercive, the
  minimization problem defined by \eqref{eq:OCP}--\eqref{eq:aDef} has a unique
  solution characterized by the
  first order optimality condition
  (cf. \cite{KS:1980:VarInequalities,ET:1999:Convex})
\begin{equation}\label{eq:OptimalityCondition1}
  \int_\O (\bar y-y_d)(y-\bar y)dx+\gamma\int_\O \bar u(u-\bar u)dx\geq0
  \qquad\forall\,(y,u)\in\bK.
\end{equation}
 \par
  Let the adjoint state $\bar p\in\HOnez$ be defined by
\begin{equation}\label{eq:AdjointState}
  a(q,\bar p)=\int_\O (\bar y-y_d)q\,dx\qquad\forall\,q\in \HOnez.
\end{equation}
 One can use \eqref{eq:PDEConstraint} and \eqref{eq:AdjointState} to write
\begin{equation}\label{eq:ASRelation}
  \int_\O (\bar y-y_d)y\,dx=
   a(y,\bar p)=\int_\O u\bar p\,dx \qquad\forall\,(y,u)\in\bK,
\end{equation}
 and then
 \eqref{eq:OptimalityCondition1}  is equivalent to the inequality
\begin{equation}\label{eq:OptimalityCondition2}
  \int_\O (\bar p+\gamma \bar u)(u-\bar u)dx\geq0 \qquad \forall\,u\in K,
\end{equation}
 where
\begin{equation*}
  K=\{u\in \LT:\,\phi_1\leq u\leq\phi_2\quad\text{a.e. in $\O$}\}.
\end{equation*}
\par
 The inequality \eqref{eq:OptimalityCondition2} is equivalent to the
 statement that $\bar u$ is the
 $\LT$-orthogonal projection of
 $-\gamma^{-1}\bar p$ on the closed convex subset $K$, i.e.,
\begin{equation}\label{eq:Truncation}
  \bar u=\max(\phi_1,\min(\phi_2,-\gamma^{-1}\bar p)),
\end{equation}
 which, in view of \eqref{eq:ConstraintRegularity}, implies in particular
  that (cf. \cite[Lemma~7.6]{GT:2001:EllipticPDE})
\begin{equation*}
  \bar u\in H^1(\O).
\end{equation*}
\par
 For the analysis of  problems with rough coefficients, it is desirable to
 keep track of the dependence of $|\bar u|_{H^1(\O)}$ on $\alpha$ and $\beta$.
 This can be achieved by using the
 estimate
\begin{equation*}
  \|\bar y-y_d\|_\LT^2\leq 2 J(y_*,u_*)
\end{equation*}
 that holds for any convenient choice of $(y_*,u_*)\in \bK$.  One can then bound
 $|\bar p|_{H^1(\O)}$ through
 \eqref{eq:AdjointState} and then obtain an estimate of $|\bar u|_{H^1(\O)}$
 through \eqref{eq:Truncation}.
\par
 For example, under the additional assumption
 $\phi_1\leq 0\leq \phi_2$
 almost everywhere in $\O$,
 we can take $(y_*,u_*)=(0,0)$ to obtain a simple bound
\begin{equation}\label{eq:AprioriBound1}
  \|\bar y-y_d\|_\LT^2\leq 2J(0,0)=\|y_d\|_\LT^2.
\end{equation}
\par
 It then follows from \eqref{eq:SpaceComparison},
 \eqref{eq:AdjointState} and
 \eqref{eq:AprioriBound1} that
\begin{equation*}
   \alpha|\bar p|_{H^1(\O)}^2\leq a(\bar p,\bar p)
                  =\int_\O (\bar y-y_d)\bar p\,dx
                  \leq \| y_d\|_\LT\|\bar p\|_\LT,
\end{equation*}
 which implies
\begin{equation}\label{eq:AprioriBound2}
  |\bar p|_{H^1(\O)}\leq (\CPF/\alpha)\|y_d\|_\LT
\end{equation}
 through the Poincar\'e-Friedrichs inequality
\begin{equation}\label{eq:PF}
  \|v\|_\LT\leq \CPF |v|_{H^1(\O)}\qquad\forall\,v\in \HOnez.
\end{equation}
\par
 Putting \eqref{eq:Truncation} and \eqref{eq:AprioriBound2} together,
 we arrive at the bound
\begin{equation*}
  |\bar u|_{H^1(\O)}\leq \max\big(|\phi_1|_{H^1(\O)},|\phi_2|_{H^1(\O)},
  \gamma^{-1}({\CPF}/\alpha)\|y_d\|_\LT\big).
\end{equation*}
\begin{remark}\label{rem:GeneralConstraints}\rm
  Under the general assumption \eqref{eq:ConstraintInequality},
  we can take $u_*=(\phi_1+\phi_2)/2$ and obtain a (more complicated)
   upper bound
  for $|u|_{H^1(\O)}$ that depends only on $\|\phi_1\|_{H^1(\O)}$,
  $\|\phi_2\|_{H^1(\O)}$, $\|y_d\|_\LT$, $\gamma^{-1}$ and $\alpha^{-1}$.
\end{remark}
\par
 Next we define
\begin{equation}\label{eq:lambda}
  \lambda=\bar p+\gamma \bar u
\end{equation}
 and obtain through \eqref{eq:Truncation} the decomposition
\begin{equation}\label{eq:LambdaDecomposition}
  \lambda=\lambda_1+\lambda_2,
\end{equation}
 where
\begin{equation*}
  \lambda_1=\max(\bar p+\gamma\phi_1,0)\in H^1(\O) \quad\text{and} \quad
  \lambda_2=\min(\bar p+\gamma\phi_2,0)\in H^1(\O)
\end{equation*}
 satisfy
\begin{subequations}\label{subeq:lambda}
\begin{align}
 \lambda_1&\geq0 ,\label{eq:lambda1Sign}\\
 \lambda_1(\phi_1-\bar u)&=0, \label{eq:lambda1Complementarity}\\
 \nabla\lambda_1&=\begin{cases}
     \nabla\bar p+\gamma\nabla \phi_1&\quad\text{in $\fA_1$}\\[2pt]
      0 &\quad\text{in $\O\setminus\fA_1$}
  \end{cases},\label{eq:lambda1Gradient}\\
  \lambda_2&\leq0,\label{eq:lambda2Sign}\\
  \lambda_2(\phi_2-\bar u)&=0, \label{eq:lambda2Complementarity}\\
  \nabla\lambda_2&=\begin{cases}
     \nabla\bar p+\gamma\nabla \phi_2&\quad\text{in $\fA_2$}\\[2pt]
     0 &\quad\text{in $\O\setminus\fA_2$}
  \end{cases}.\label{eq:lambda2Gradient}
\end{align}
\end{subequations}
 Here the active set $\fA_j$ is the closure in $\O$ of the set of
 the Lebesgue points  where $\bar u-\phi_j=0$.
%
%%%%%%%%%%%%%%%%%%%%%%%%%%%%%
\section{A Discretization of the Optimal Control Problem}\label{sec:FEMs}
 Let $\cT_\rho$ be a simplicial/quadrilateral triangulation of $\O$ with
 mesh size $\rho$ and $W_\rho\subset \LT$
 be the space of piecewise constant functions with
 respect to $\cT_\rho$.  The optimal control $\bar u$ will be approximated
 by functions in $W_\rho$, while
 the approximation of $\bar y$ comes from a subspace $V_*$ of $\HOnez$.
\begin{remark}\label{rem:General}\rm
  By allowing $V_*$ to be an arbitrary subspace of $H^1_0(\O)$, the analysis
  developed below can be applied to
  standard finite element methods and multiscale finite element methods.
\end{remark}
\par
 The discrete problem is to find
 \begin{equation}\label{eq:DOCP}
   (\bys,\bur)=\argmin_{(y_*,u_\rho)\in \bKsr} J(y_*,u_\rho),
 \end{equation}
  where the closed convex subset $\bKsr$ of $V_*\times W_\rho$ is
  defined by  the following conditions:
\begin{alignat}{3}
   a(y_*,z_*)&=\int_\O u_\rho z_* dx&\qquad&\forall\,z_*\in V_*,
   \label{eq:DPDE}\\
    Q_\rho \phi_1\leq\,&u_\rho\leq Q_\rho\phi_2&\qquad&\text{a.e. in $\O$},
    \label{eq:DCC}
\end{alignat}
 and $Q_\rho$ is the orthogonal projection from $\LT$ onto $W_\rho$.
\par
 We have a standard interpolation error estimate
  (cf. \cite{Ciarlet:1978:FEM,BScott:2008:FEM})
\begin{equation}\label{eq:QEstimate}
  \|\zeta-Q_\rho\zeta\|_\LT\leq C_\maltese\rho|\zeta|_{H^1(\O)},
\end{equation}
 where the positive constant $C_\maltese$ only depends on the shape regularity
  of $\cT_\rho$.
\par
 Since $Q_\rho u$ satisfies \eqref{eq:DCC} for any $u$ that satisfies
 \eqref{eq:ControlConstraints}, the set $\bK_{*,\rho}$ is
 nonempty and the discrete convex minimization problem
 has a unique solution characterized by the first order optimality condition
\begin{equation}\label{eq:DOC}
  \int_\O (\bys-y_d)(y_*-\bys)dx+
  \gamma\int_\O \bur(u_\rho-\bur)dx\geq0 \qquad
  \forall\, (y_*,u_\rho)\in \bKsr.
\end{equation}
\par
 The error analysis for $(\bys,\bur)$ was carried out in the pioneering
 work \cite{Falk:1973:Control} on finite element methods for elliptic
 optimal control problems.
 Here we present a self-contained treatment that is suitable for the
 analysis of the multiscale
 finite element method in Section~\ref{sec:DDLOD}.
\par
 The following lemma is useful for the error analysis.
\begin{lemma}\label{lem:EasyEstimate}
  Let $g\in\LT$ and $w_*\in V_*$ satisfy
\begin{equation*}
  a(w_*,v_*)=\int_\O g v_* dx\qquad\forall\,v_*\in V_*.
\end{equation*}
 Then we have
\begin{align}
  \|w_*\|_\LT&\leq (\CPF^2/\alpha)\|g\|_\LT,\label{eq:EasyEstimate}\\
  \|w_*\|_a&\leq (\CPF/\sqrt\alpha)\|g\|_\LT.\label{eq:EasyEnergyEstimate}
\end{align}
\end{lemma}
\begin{proof}
  It follows from \eqref{eq:SpaceComparison}, \eqref{eq:PF} and the
  Cauchy-Schwarz inequality that
\begin{align*}
  \|w_*\|_\LT^2\leq \CPF^2|w_*|_{H^1(\O)}^2
               &\leq (\CPF^2/\alpha)a(w_*,w_*)\\
  &= (\CPF^2/\alpha)\int_\O gw_*dx\leq (\CPF^2/\alpha)\|g\|_\LT\|w_*\|_\LT,
\end{align*}
 which implies \eqref{eq:EasyEstimate}.
\par
 The estimate \eqref{eq:EasyEnergyEstimate} also follows from
 \eqref{eq:SpaceComparison},
 \eqref{eq:PF} and the Cauchy-Schwarz inequality:
\begin{equation*}
  \|w_*\|_a^2=a(w_*,w_*)=\int_\O gw_*dx\leq \|g\|_\LT \|w_*\|_\LT\leq
  \|g\|_\LT(\CPF/\sqrt\alpha)\|w_*\|_a.
\end{equation*}

\end{proof}
\par
 We will include the approximation of $\bar p$ by $\bps$ in the error analysis, where
 $\bps\in V_*$ is defined by
\begin{equation}\label{eq:bps}
  a(q_*,\bps)=\int_\O (\bys-y_d)q_* dx\qquad\forall\,q_*\in V_*.
\end{equation}
\begin{theorem}\label{thm:AbstractErrorEstimate}
  There exists a positive constant $C_\dag$, depending only on
  $\|y_d\|_\LT$, $\|\phi_1\|_{H^1(\O)}$, $\|\phi_2\|_{H^1(\O)}$,
  $\gamma^{-1}$, $\alpha^{-1}$
 and the shape regularity of $\cT_\rho$, such that
\begin{equation}\label{eq:AbstractErrorEstimate}
  \|\bar y-\bys\|_\LT+\|\bar u-\bur\|_\LT+\|\bar p-\bps\|_\LT
  \leq C_\dag(\|\bar y-\dys\|_\LT+
  \|\bar p-\dps\|_\LT+\rho),
\end{equation}
 where $\dys,\dps\in V_*$ are defined by
\begin{alignat}{3}
  a(\dys,z_*)&=\int_\O \bar u z_* dx&\qquad&\forall\,z_*\in V_*,\label{eq:dys}\\
  a(q_*,\dps)&=\int_\O (\bar y-y_d)q_* dx&\qquad&\forall\,q_*\in V_*.\label{eq:dps}
\end{alignat}
\end{theorem}
\begin{proof}  First we note the following analog of \eqref{eq:ASRelation}:
\begin{equation}\label{eq:DiscreteASRelation}
  \int_\O (\bar y-y_d)y_*dx=
   a(y_*,\dps)=\int_\O u_\rho\dps dx \qquad\forall\,(y_*,u_\rho)\in\bKsr
\end{equation}
 by \eqref{eq:DPDE} and \eqref{eq:dps}.
\par
 Let $(\tys,\tur)\in\bKsr$ be defined by
\begin{equation}\label{eq:tur}
 \tur=Q_\rho\bar u
\end{equation}
 and
\begin{equation}\label{eq:tys}
  a(\tys,z_*)=\int_\O \tur z_*\,dx \qquad\forall\,z_*\in V_*.
\end{equation}
\par
 It follows from \eqref{eq:QEstimate} and \eqref{eq:tur} that
\begin{equation}\label{eq:turEstimate}
  \|\bar u-\tur\|_\LT\leq C_\maltese\rho|\bar u|_{H^1(\O)}.
\end{equation}
\par
 We have
\begin{align}\label{eq:SError1}
  &\|\bar y-\bys\|_\LT^2+\gamma\|\bar u-\bur\|_\LT^2\notag\\
  &\hspace{30pt}=\int_\O (\bar y-\bys)(\bar y-\tys)dx+
  \gamma\int_\O (\bar u-\bur)(\bar u-\tur)dx\\
  &\hspace{60pt}+
   \int_\O (\bar y-\bys)(\tys-\bys)dx+
   \gamma\int_\O (\bar u-\bur)(\tur-\bur)dx,\notag
\end{align}
 and, in view of \eqref{eq:lambda}, \eqref{eq:DOC} and
 \eqref{eq:DiscreteASRelation},
\begin{align}\label{eq:SError2}
  &\int_\O (\bar y-\bys)(\tys-\bys)dx+
  \gamma\int_\O (\bar u-\bur)(\tur-\bur)dx\notag\\
  &\hspace{30pt}=\int_\O \bar y(\tys-\bys)dx+
  \gamma\int_\O \bar u(\tur-\bur)dx\notag\\
  &\hspace{60pt}-\int_\O \bys(\tys-\bys)dx-
  \gamma\int_\O \bur(\tur-\bur)dx\\
   &\hspace{30pt}\leq\int_\O (\bar y-y_d)(\tys-\bys)dx
   +\gamma\int_\O \bar u(\tur-\bur)dx\notag\\
   &\hspace{30pt}=\int_\O (\dps+\gamma\bar u)(\tur-\bur)dx\notag\\
   &\hspace{30pt}=\int_\O \lambda(\tur-\bur)dx
   +\int_\O (\dps-\bar p)(\tur-\bur)dx.\notag
\end{align}
\par
 We can bound the first integral on the right-hand side
 of \eqref{eq:SError2} by \eqref{eq:ConstraintRegularity},
 Remark~\ref{rem:GeneralConstraints},
 \eqref{eq:LambdaDecomposition}, \eqref{subeq:lambda}, \eqref{eq:DCC},
 \eqref{eq:QEstimate} and \eqref{eq:turEstimate}:
\begin{align}\label{eq:SError3}
  \int_\O \lambda(\tur-\bur)dx&=\int_\O\lambda_1(\tur-\bur)dx
  +\int_\O\lambda_2(\tur-\bur)dx\notag\\
     &=\int_\O \lambda_1(Q_\rho\bar u-\bar u)dx
     +\int_\O \lambda_2(Q_\rho\bar u-\bar u)dx\notag\\
     &\hspace{40pt}+\int_\O \lambda_1(\bar u-\phi_1)dx
     +\int_\O \lambda_2(\bar u-\phi_2)dx\notag\\
     &\hspace{60pt}+\int_\O \lambda_1(\phi_1-Q_\rho\phi_1)dx
     +\int_\O \lambda_2(\phi_2-Q_\rho\phi_2)dx\notag\\
     &\hspace{80pt}+\int_\O \lambda_1(Q_\rho\phi_1-\bur)dx+
     \int_\O \lambda_2(Q_\rho\phi_2-\bur)dx\notag\\
     &\leq \int_\O\lambda_1 (Q_\rho\bar u-\bar u)dx+
     \int_\O\lambda_2(Q_\rho\bar u-\bar u)dx\\
      &\hspace{40pt}+\int_\O\lambda_1(\phi_1-Q_\rho\phi_1)dx+
      \int_\O\lambda_2(\phi_2-Q_\rho\phi_2)dx\notag\\
      &=\int_\O(\lambda_1-Q_\rho\lambda_1) (Q_\rho\bar u-\bar u)dx+
      \int_\O(\lambda_2-Q_\rho\lambda_2)(Q_\rho\bar u-\bar u)dx\notag\\
 &\hspace{40pt}+\int_\O(\lambda_1-Q_\rho\lambda_1)(\phi_1-Q_\rho\phi_1)dx+
     \int_\O(\lambda_2-Q_\rho\lambda_2)(\phi_2-Q_\rho\phi_2)dx\notag\\
      &\leq C_1\rho^2.\notag
\end{align}
\par
 For the second integral on the right-hand side of \eqref{eq:SError2},
  we have
\begin{align}\label{eq:SError4}
  \int_\O (\dps-\bar p)(\tur-\bur)dx&
  \leq \|\bar p-\dps\|_\LT(\|\tur-\bar u\|_\LT+\|\bar u-\bur\|_\LT)\\
  &\leq \|\bar p-\dps\|_\LT\big(C_\maltese\rho|\bar u|_{H^1(\O)}
  +\|\bar u-\bur\|_\LT\big)\notag
\end{align}
 by  \eqref{eq:turEstimate}.
\par
 It follows from  Remark~\ref{rem:GeneralConstraints},
 \eqref{eq:turEstimate},  and
 \eqref{eq:SError1}--\eqref{eq:SError4} that
\begin{align*}
  &\|\bar y-\bys\|_\LT^2+\gamma\|\bar u-\bur\|_\LT^2\\
  &\hspace{70pt}\leq \|\bar y-\bys\|_\LT\|\bar y-\tys\|_\LT
 +\gamma\|\bar u-\bur\|_\LT C_\maltese\rho|\bar u|_{H^1(\O)}\\
   &\hspace{120pt}+C_1\rho^2+C_2 \rho\|\bar p-\dps\|_\LT
   +\|\bar p-\dps\|_\LT\|\bar u-\bur\|_\LT,
\end{align*}
 which together with
 the inequality of arithmetic and geometric means implies
\begin{equation}\label{eq:SError5}
  \|\bar y-\bys\|_\LT+\|\bar u-\bur\|_\LT\leq C_3\big(\|\bar y-\tys\|_\LT
  +\|\bar p-\dps\|_\LT+\rho).
\end{equation}
\par
 On the other hand, we have
\begin{equation*}
  a(\dys-\tys,z_*)=\int_\O (\bar u-\tur)z_*dx\qquad\forall\,z_h\in V_*,
\end{equation*}
 by \eqref{eq:dys} and \eqref{eq:tys}, and hence
\begin{equation}\label{eq:SError6}
  \|\dys-\tys\|_\LT\leq (\CPF^2/\alpha)\|\bar u-\tur\|_\LT
  \leq (\CPF^2/\alpha)C_\maltese\rho|\bar u|_{H^1(\O)}
\end{equation}
 by \eqref{eq:QEstimate} and Lemma~\ref{lem:EasyEstimate}.
\par
 Putting \eqref{eq:SError5} and \eqref{eq:SError6} together, we arrive at
 the estimate
\begin{equation}\label{eq:SError7}
  \|\bar y-\bys\|_\LT+\|\bar u-\bur\|_\LT\leq C_4\big(\|\bar y-\dys\|_\LT
  +\|\bar p-\dps\|_\LT+\rho).
\end{equation}

\par
 For the estimate of $\bar p-\bps$, we begin with
\begin{equation}\label{eq:SError8}
  \|\bar p-\bps\|_\LT\leq \|\bar p-\dps\|_\LT+\|\dps-\bps\|_\LT
\end{equation}
 and note that
\begin{equation}\label{eq:SError10}
  a(q_*,\dps-\bps)=\int_\O (\bar y-\bys)q_*dx \qquad\forall\,q_*\in V_*
\end{equation}
 by \eqref{eq:bps} and \eqref{eq:dps}, which implies
\begin{equation}\label{eq:SError9}
  \|\dps-\bps\|_\LT\leq (\CPF^2/\alpha)\|\bar y-\bys\|_\LT
\end{equation}
 through Lemma~\ref{lem:EasyEstimate}.
\par
 The estimate \eqref{eq:AbstractErrorEstimate} follows from
 \eqref{eq:SError7}--\eqref{eq:SError9}.
\end{proof}
\par
 The following result shows that the estimate
 \eqref{eq:AbstractErrorEstimate} is a tight estimate.
\begin{theorem}\label{thm:ReverseErrorEstimate}
  There exists a positive constant $C_\ddag$, depending only on
  $\alpha^{-1}$, such that
\begin{equation}\label{eq:ReverseErrorEstiamte}
 \|\bar y-\dys\|_\LT+\|\bar p-\dps\|_\LT\leq C_\ddag\big(\|\bar y-\bys\|_\LT
  +\|\bar u-\bur\|_\LT+\|\bar p-\bps\|_\LT\big),
\end{equation}
 where $\dys$ $($resp., $\dps)$ is defined by \eqref{eq:dys}
  $($resp., $\eqref{eq:dps})$.
\end{theorem}
\begin{proof} We have
\begin{equation}\label{eq:REstimate1}
  \|\bar y-\dys\|_\LT\leq \|\bar y-\bys\|_\LT+\|\bys-\dys\|_\LT,
\end{equation}
 and
\begin{equation*}
  a(\bys-\dys,z_*)=\int_\O (\bur-\bar u)z_* dx\qquad\forall\,z_*\in V_*
\end{equation*}
 by \eqref{eq:DPDE} and \eqref{eq:dys}, which implies
\begin{equation}\label{eq:REstiamte2}
  \|\bys-\dys\|_\LT\leq (\CPF^2/\alpha)\|\bar u-\bur\|_\LT
\end{equation}
 by Lemma~\ref{lem:EasyEstimate}.
 \par
  Similarly we have
 \begin{align}\label{eq:REstimate3}
   \|\bar p-\dps\|_\LT&\leq \|\bar p-\bps\|_\LT+\|\bps-\dps\|_\LT\\
       &\leq \|\bar p-\bps\|_\LT+(\CPF^2/\alpha)\|\bar y-\bys\|_\LT\notag
 \end{align}
 by \eqref{eq:SError9}.
 \par
 The estimate \eqref{eq:ReverseErrorEstiamte} follows from
 \eqref{eq:REstimate1}--\eqref{eq:REstimate3}.
\end{proof}
\par
 It is straightforward to derive error estimates in the energy norm
 from  the estimate \eqref{eq:AbstractErrorEstimate}.
\begin{theorem}\label{thm:AbstractEnergyError}
  There exists a positive constant $C_\S$, depending only on
  $\|y_d\|_\LT$, $\|\phi_1\|_{H^1(\O)}$, $\|\phi_2\|_{H^1(\O)}$,
  $\gamma^{-1}$, $\alpha^{-1}$
 and the shape regularity of $\cT_\rho$, such that
\begin{equation}\label{eq:AbastractEnergyError}
  \|\bar y-\bys\|_a+\|\bar p-\bps\|_a\leq C_\S \big(\|\bar y-\dys\|_a
  +\|\bar p-\dps\|_a+\rho\big),
\end{equation}
 where $\dys, \dps\in V_*$ are defined in \eqref{eq:dys} and \eqref{eq:dps}.
\end{theorem}
\begin{proof}
  We have
\begin{equation}\label{eq:AEnergy1}
  \|\bar y-\bys\|_a+\|\bar p-\bps\|_a\leq \|\bar y-\dys\|_a
  +\|\bar p-\dps\|_a+\|\dys-\bys\|_a
  +\|\dps-\bps\|_a.
\end{equation}
\par
 It follows from \eqref{eq:DPDE} and \eqref{eq:dys} that
\begin{equation*}
  a(\dys-\bys,z_*)=\int_\O (\bar u-\bur)z_*dx\qquad\forall\,z_*\in V_*,
\end{equation*}
 and hence
\begin{equation}\label{eq:AEnergy2}
  \|\dys-\bys\|_a\leq (\CPF/\sqrt\alpha)\|\bar u-\bur\|_\LT
\end{equation}
 by Lemma~\ref{lem:EasyEstimate}.
\par
 Similarly the relation \eqref{eq:SError10} and Lemma~\ref{lem:EasyEstimate}
 imply
\begin{equation}\label{eq:AEnergy3}
  \|\dps-\bps\|_a\leq (\CPF/\sqrt\alpha)\|\bar y-\bys\|_\LT.
\end{equation}
\par
 The estimate \eqref{eq:AbastractEnergyError} is obtained by
 combining \eqref{eq:AbstractErrorEstimate},
 \eqref{eq:AEnergy1}--\eqref{eq:AEnergy3} and the relation
\begin{equation*}
  \|\bar y-\dys\|_\LT+\|\bar p-\dps\|_\LT\leq(\CPF/\sqrt\alpha)\big(
   \|\bar y-\dys\|_a+\|\bar p-\dps\|_a\big)
\end{equation*}
 that follows from \eqref{eq:SpaceComparison} and \eqref{eq:PF}.
\end{proof}
\begin{remark}\label{rem:GalerkinProjection}\rm
  Note that \eqref{eq:AdjointState} and \eqref{eq:dps} imply
  $\dps\in V_*$ is the projection of
  $\bar p$ with respect to the bilinear form $a(\cdot,\cdot)$.
   Therefore we have
\begin{equation*}
   \|\bar p-\dps\|_a=\inf_{q_*\in V_*}\|\bar p-q_*\|_a.
\end{equation*}
 Similarly we have
\begin{equation*}
  \|\bar y-\dys\|_a=\inf_{z_*\in V_*}\|\bar y-z_*\|_a
\end{equation*}
 by \eqref{eq:PDEConstraint} and \eqref{eq:dys}.
\end{remark}
\par
  Let $V_*=V_h$ be the $P_1/Q_1$ finite element space associated
 with a simplicial/quadrilateral triangulation $\cT_h$ of $\O$ with mesh size $h$
   and let $(\bys,\bur,\bps)$ be written as $(\byh,\buh,\bph)$. The estimate
   \eqref{eq:AbstractErrorEstimate}
   becomes
\begin{equation}\label{eq:StandardErrorEstimate}
  \|\bar y-\byh\|_\LT+\|\bar u-\buh\|_\LT+\|\bar p-\bph\|_\LT
  \leq C_\dag(\|\bar y-\dyh\|_\LT+
  \|\bar p-\dph\|_\LT+\rho),
\end{equation}
 where $\dyh,\dph\in V_h$ are defined by
\begin{alignat}{3}
  a(\dyh,z_h)&=\int_\O \bar u z_h dx&\qquad&\forall\,z_h\in V_h,\label{eq:dyh}\\
  a(q_h,\dph)&=\int_\O (\bar y-y_d)q_h dx&\qquad&\forall\,q_h\in V_h,\label{eq:dph}
\end{alignat}
 and the estimate \eqref{eq:AbastractEnergyError} becomes
\begin{equation}\label{eq:StandardEnergyError}
  \|\bar y-\byh\|_a+\|\bar p-\bph\|_a\leq C_\S\big(\|\bar y-\dyh\|_a+\|\bar p-\dph\|_a+\rho\big).
\end{equation}
\par
 In the case where $\cA$ is the identity matrix and $\O$ is
 convex, we have $\bar y,\bar p \in H^2(\O)$ by the elliptic regularity theory for
 polygonal domains (cf. \cite{Grisvard:1985:EPN,Dauge:1988:EBV,MR:2010:Polyhedral}).
 It follows from Remark~\ref{rem:GalerkinProjection},
 \eqref{eq:StandardErrorEstimate}  and a
 standard duality argument (cf. \cite{Ciarlet:1978:FEM,BScott:2008:FEM}) that
\begin{equation}\label{eq:Smooth}
  \|\bar y-\byh\|_\LT+\|\bar u-\buh\|_\LT+\|\bar p-\bph\|_\LT\leq C(h^2+\rho).
\end{equation}
 In this case the estimate \eqref{eq:StandardEnergyError} yields
\begin{equation}\label{eq:SmoothEnergy}
  |\bar y-\byh|_{H^1(\O)}+|\bar p-\bph|_{H^1(\O)}\leq C(h+\rho).
\end{equation}
\par
 In the case of rough coefficients, we can derive from \eqref{eq:PF},
 \eqref{eq:AbstractErrorEstimate}, \eqref{eq:AbastractEnergyError} and
  Remark~\ref{rem:GalerkinProjection} that
\begin{align*}
  &\|\bar y-\byh\|_\LT+\|\bar u-\buh\|_\LT+\|\bar p-\bph\|_\LT
    +\|\bar y-\byh\|_a+\|\bar p-\bph\|_a\\
  &\hspace{80pt}\leq C\Big(\inf_{z_h\in V_h}\|\bar y-z_h\|_a
  +\inf_{q_h\in V_h}\|\bar p-q_h\|_a+\rho\Big),
\end{align*}
 which implies
 $$\lim_{h,\rho\downarrow0}\big(\|\bar y-\byh\|_\LT+\|\bar u-\buh\|_\LT
 +\|\bar p-\bph\|_\LT
 +\|\bar y-\byh\|_a+\|\bar p-\bph\|_a\big)=0.$$
 However the convergence with respect to $h$ can be very slow.  Therefore
  a satisfactory
 approximate solution of the
 optimal control problem obtained by standard finite element methods will
 require a very fine mesh $\cT_h$.
\par
 Below we will show that it is possible to recover on coarse meshes
  a performance similar
 to \eqref{eq:Smooth} and \eqref{eq:SmoothEnergy} for
 rough coefficients and general
 $\O$ provided that one takes a multiscale finite element space to be $V_*$.
%%%%%%%%%%%%%%%%%%%%%%%%%%%%%%%%%%%%%%%%%%%%%%%%
\section{A DD-LOD Multiscale Finite Element Method}\label{sec:DDLOD}
 First we recall the construction of the multiscale finite element space
 from \cite{BGS:2021:LOD}.
 It begins with
 a simplicial/quadrilateral triangulation  $\cT_H$ of $\O$, and  a
  refinement $\cT_h$ ($h\ll H$) of
 $\cT_H$.  The $P_1/Q_1$ finite element subspace of $H^1_0(\O)$ associated with
 $\cT_H$ (resp., $\cT_h$) is denoted by $V_H$ (resp., $V_h$).
\par
 The first step is to construct a projection operator
 $\PiH:H^1_0(\O)\longrightarrow V_H$ such that
\begin{equation*}
  \frac{1}{H}\|v-\PiH v\|_\LT+|\PiH v|_{H^1(\O)}\leq C_\flat|v|_{H^1(\O)} \qquad
  \forall\,v\in H^1_0(\O).
\end{equation*}
\begin{remark}\label{rem:PiH}\rm
  The operator $\PiH$ in \cite{BGS:2021:LOD} is constructed by taking the averages of
  local $L_2$ projections.
  There are other constructions that are adapted to the coefficient matrix $\cA(x)$
  (cf. \cite{PS:2016:Contrast,HM:2017:Contrast}).
\end{remark}
\par
 Let $\KER=\{v\in V_h:\PiH v=0\}$ be the kernel of $\PiH$ in $V_h$ and
 the correction operator $\Ch:V_h\longrightarrow\KER$ be the projection
 operator with respect to
 $a(\cdot,\cdot)$, i.e.,
\begin{equation*}%\label{eq:Correction}
  a(\Ch v,w)=a(v,w) \qquad\forall\,w\in \KER.
\end{equation*}
\par
 The multiscale finite element space $\MS\subset V_h$ is the orthogonal complement of
 $\KER$ with respect to
 $a(\cdot,\cdot)$.  Let $\phi_1,\ldots,\phi_m$ be the standard nodal basis functions
 of $V_H$ associated with the interior
 vertices $p_1,\ldots,p_m$ of $\cT_H$.  Then $\MS$ is spanned by
 $\phi_1-\Ch\phi_1,\ldots,\phi_m-\Ch\phi_m$.
 The performance of the finite element method based on $\MS$ for the problem
\begin{equation}\label{eq:RBVP}
  a(u,v)=\int_\O fv\,dx \qquad\forall\,v\in \HOnez
\end{equation}
 with rough
 coefficients
 is similar to the performance of $V_H$ for problems with smooth coefficients on
 convex domains
  (cf. \cite{MP:2014:LOD,MP:2021:LOD}).
 However, the construction of $\MS$ requires solving $m$ problems on the fine mesh
 $\cT_h$, which is expensive.
\par
 The localized orthogonal decomposition (LOD) method is
 based on replacing the correction $\Ch\phi_i$ by a correction computed in a subdomain
 consisting of
 a certain number of layers of elements from $\cT_H$ around $p_i$.  It significantly
 reduces the computational cost and
 at the same time it preserves the good approximation property of $\MS$
 because the function
 $\Ch\phi_i$ decays exponentially
 away from $p_i$ (cf. \cite{MP:2014:LOD,MP:2021:LOD,AHP:2021:Acta}).
\par
 The multiscale finite element method from \cite{BGS:2021:LOD} is a variant
 of the LOD method which is
 based on the ideas in
 \cite{KPY:2018:LOD}.  It computes an approximate solution $\Chk\phi_i$ of
 the corrector equation
\begin{equation*}
  a(\Ch\phi_i,w)=a(\phi_i,w) \qquad \forall\,w \in \KER
\end{equation*}
 by applying $k$ iterations of a preconditioned conjugate gradient (PCG) method
 with initial guess $0$.
  The theory of PCG (cf. \cite{Saad:2003:IM}) implies that the convergence of
  $\Chk\phi_i$ to
  $\Ch\phi_i$ in $\|\cdot\|_a$ is approximately
 $q^k$, where $q\in (0,1)$ depends on the condition number of the preconditioned
 operator.
\par
 The key is to use an additive Schwarz domain decomposition preconditioner
 (cf. \cite{TW:2005:DD}) where the subdomains
 are small patches $\omega_i$ around $p_i$ so that $\Chk\phi_i$ is supported on
 a subdomain obtained by
 adding approximately $2k$ layers of elements from $\cT_H$  around $\omega_i$, i.e.,
 $\Chk\phi_i$ is also a localized
 correction of $\phi_i$.
 The computation of $\Chk\phi_i$ only involves solving local small problems
  and $\|\Ch\phi_i-\Chk\phi_i\|_a=O(H)$ provided $k$ is proportional to $|\ln H|$.
\par
 The multiscale finite element space
 $\MSV\subset V_h$ is spanned by $\phi_1-\Chk\phi_1,\ldots,\phi_m-\Chk\phi_m$.
 We will refer to
 it as the DD-LOD multiscale finite element space.  The corresponding finite
 element method
 for \eqref{eq:RBVP} can be viewed as a
 reduced order method, where
 the functions $\Chk\phi_1,\ldots,\Chk\phi_m$ are computed off-line.
 The on-line computation
 only involves solving an $m\times m$ system.
\par
 The following is the main result from \cite{BGS:2021:LOD} whose
 derivation only
 involves basic results
 from finite element methods, domain decomposition methods and numerical
 linear algebra.
\begin{lemma}\label{lem:ETNA}
  Let $f\in\LT$, $y_h\in V_h$ and $\MSLy\in\MSV$ such that
\begin{alignat*}{3}
  a(y_h,z_h)&=\int_\O fz_h dx&\qquad&\forall\,z_h\in V_h, \\
  a(\MSLy,\MSLz)&=\int_\O f\MSLz dx &\qquad&\forall\,\MSLz\in \MSV.
\end{alignat*}
 There exists a positive constant $C_\sharp$ depending on the shape regularity of
   $\cT_H$ but independent of $\alpha$, $\beta$, $h$ and $H$, such that
\begin{align*}
  \|y_h-\MSLy\|_a&\leq (C_\sharp/\sqrt{\alpha})H\|f\|_\LT,\\
  \|y_h-\MSLy\|_\LT&\leq (C_\sharp/\sqrt{\alpha})^2H^2\|f\|_\LT,
\end{align*}
 provided $k=\lceil -j \ln H\rceil$ for a sufficiently large $j$.
\end{lemma}
\begin{remark}\label{rem:j}\rm
  The magnitude of $j$ depends on the condition number of the
   preconditioned operator
  in the PCG algorithm.
\end{remark}
%%%%%%%%%%%%%%%%%%%%%%%%%%%%%%
\par
 The DD-LOD finite element method for \eqref{eq:OCP}--\eqref{eq:aDef}
 is defined by
 \eqref{eq:DOCP}--\eqref{eq:DCC},
 where $V_*=\MSV$
 and its solution is denoted
 by $(\bMSLy,\bar u_\rho)$.
\par
 We also include the approximation of $\bar p$ by $\bMSLp$ in the error
 analysis of the multiscale finite element method, where
 $\bMSLp\in\MSV$ is defined by
\begin{equation}\label{eq:bMSLp}
  a(\MSLq,\bMSLp)=\int_\O (\bMSLy-y_d)\MSLq \,dx\qquad\forall\,\MSLq\in\MSV.
\end{equation}
\begin{remark}\label{rem:Notation}\rm
  Strictly speaking $\bMSLy$ and $\bMSLp$ also depend on $\rho$ and
  $\bar u_\rho$ also depends on $h$, $H$ and $k$.  These dependencies
   are suppressed
  for the sake of readability.
\end{remark}
\begin{theorem}\label{thm:DDLODError}
  There exists a positive constant $C_\natural$, depending
  only on
  $\|y_d\|_\LT$, $\|\phi_1\|_{H^1(\O)}$, $\|\phi_2\|_{H^1(\O)}$,
  $\gamma^{-1}$, $\alpha^{-1}$
 and the shape regularities of $\cT_H$ and $\cT_\rho$, such that
\begin{align}\label{eq:DDLODError}
    &\|\bar y-\bMSLy\|_\LT+\|\bar u-\bar u_\rho\|_\LT+\|\bar p-\bMSLp\|_\LT\\
    &\hspace{80pt}\leq
  C_\natural\big(\|\bar y-\byh\|_\LT+\|\bar u-\buh\|_\LT+\|\bar p-\bph\|_\LT
  + H^2 + \rho \big),\notag
\end{align}
 where $(\byh,\buh,\bph)$ is the approximation of $(\bar y,\bar u,\bar p)$
 obtained by using
 the standard finite element space $V_h\times W_\rho$ in the discretization
 defined by \eqref{eq:DOCP}--\eqref{eq:DCC}.
\end{theorem}
\begin{proof}
\par
 We apply Theorem~\ref{thm:AbstractErrorEstimate} (with $V_*=\MSV$)
 to obtain
\begin{align}\label{eq:Analog1}
  &\|\bar y-\bMSLy\|_\LT+\|\bar u-\bar u_\rho\|_\LT+\|\bar p-\bMSLp\|_\LT\\
  &\hspace{50pt}\leq
  C_\dag\big(\|\bar y-\dMSLy\|_\LT+\|\bar p-\dMSLp\|_\LT+ \rho\big),\notag
\end{align}
 where $\dMSLy\in \MSV$ (resp., $\dMSLp\in \MSV$) is the analog of
 $\dys$ in \eqref{eq:dys}
  (resp., $\dps$ in \eqref{eq:dps}), i.e.,
 $\dMSLy$ is defined by
\begin{equation}\label{eq:dMSLy}
a(\dMSLy,\MSLz)=\int_\O \bar u \MSLz dx\qquad\forall\,\MSLz\in \MSV,
\end{equation}
  and $\dMSLp$ is defined by
\begin{equation}\label{eq:dMSLp}
 a(\MSLq,\dMSLp)=\int_\O(\bar y-y_d)\MSLq dx\qquad\forall\,\MSLq\in \MSV.
\end{equation}
\par
 Let  $\dyh\in V_h$ $($resp., $\dph\in V_h)$ be defined by
  \eqref{eq:dyh} (resp., \eqref{eq:dph}).  According to
  Theorem~\ref{thm:ReverseErrorEstimate}, we have
\begin{equation}\label{eq:LODDDEst0}
  \|\bar y-\dyh\|_\LT+\|\bar p-\dph\|_\LT
  \leq C_\ddag\big(\|\bar y-\byh\|_\LT+\|\bar u-\buh\|_\LT+
    \|\bar p-\bph\|_\LT\big).
\end{equation}
\par
 On the other hand,
 in view of Lemma~\ref{lem:ETNA}, we have
\begin{equation}\label{eq:LODDDEst1}
  \|\dyh-\dMSLy\|_\LT\leq (C_\sharp/\sqrt{\alpha})^2H^2\|\bar u\|_\LT
\end{equation}
 by \eqref{eq:dyh} and \eqref{eq:dMSLy}, and
\begin{equation}\label{eq:LODDDEst2}
  \|\dph-\dMSLp\|_\LT\leq (C_\sharp/\sqrt{\alpha})^2H^2\|\bar y-y_d\|_\LT
\end{equation}
 by \eqref{eq:dph} and \eqref{eq:dMSLp}.
\par
 The estimate \eqref{eq:DDLODError} follows from \eqref{eq:Analog1},
 \eqref{eq:LODDDEst0}--\eqref{eq:LODDDEst2} and the triangle inequality.
\end{proof}
\begin{remark}\label{rem:Comparison}\rm
 The estimate \eqref{eq:DDLODError}
 indicates that
 up to an $O(H^2+\rho)$ error the approximation of $(\bar y,\bar u,\bar p)$ by
 $(\bMSLy,\bar u_\rho,\bMSLp)$ is as good as the approximation by
 $(\byh,\buh,\bph)$.  On the other hand, by comparing \eqref{eq:Smooth} and
 \eqref{eq:DDLODError}, we
  can also say that, up to the fine scale error,
   the performance of the multiscale finite element method on coarse
   meshes (with respect to the $\LT$ norm) is similar to
  the performance of standard finite element methods for problems with
  smooth coefficients on
  convex domains.
\end{remark}
\par
 We also have error estimates in the energy norm.
\begin{theorem}\label{thm:EnergyErrors}
  There exists a positive constant $C_\diamond$, depending
  only on
  $\|y_d\|_\LT$, $\|\phi_1\|_{H^1(\O)}$, $\|\phi_2\|_{H^1(\O)}$,
  $\gamma^{-1}$, $\alpha^{-1}$
 and the shape regularities of $\cT_H$ and $\cT_\rho$, such that
\begin{equation}\label{eq:EnergyErrors}
  \|\bar y-\bMSLy\|_a+\|\bar p-\bMSLp\|_a\leq C_\diamond
  \big(\|\bar y-\byh\|_a+\|\bar p-\bph\|_a + H + \rho\big),
\end{equation}
 where $(\byh,\bph)$ is the approximation of $(\bar y,\bar p)$ obtained by using
 the standard finite element space $V_h\times W_\rho$ in the discretization
 defined by \eqref{eq:DOCP}--\eqref{eq:DCC}.
\end{theorem}
\begin{proof} It follows from Theorem~\ref{thm:AbstractEnergyError} that
\begin{equation}\label{eq:Energy1}
  \|\bar y-\bMSLy\|_a+\|\bar p-\bMSLp\|_a
  \leq C_\S\big(\|\bar y-\dMSLy\|_a+\|\bar p-\dMSLp\|_a+\rho\big),
\end{equation}
 where $\dMSLy,\dMSLp\in \MSV$ are defined by \eqref{eq:dMSLy} and \eqref{eq:dMSLp}.
\par
 Let  $\dyh\in V_h$ $($resp., $\dph\in V_h)$ be defined by
  \eqref{eq:dyh} (resp., \eqref{eq:dph}).   In view of Lemma~\ref{lem:ETNA},
  we have
\begin{equation}\label{eq:Energy2}
  \|\dyh-\dMSLy\|_a\leq \big(C_\sharp/\sqrt{\alpha})H\|\bar u\|_\LT
\end{equation}
   by \eqref{eq:dyh} and  \eqref{eq:dMSLy}, and also
\begin{equation}\label{eq:Energy3}
\|\dph-\dMSLp\|_a \leq C_\sharp/\sqrt{\alpha})H\|\bar y-y_d\|_\LT
\end{equation}
 by \eqref{eq:dph} and \eqref{eq:dMSLp}.
\par
 Finally we note that
\begin{align}\label{eq:Energy4}
  \|\bar y-\dyh\|_a\leq \|\bar y-\byh\|_a \quad\text{and}\quad
  \|\bar p-\dph\|_a\leq \|\bar p-\bph\|_a
\end{align}
 by Remark~\ref{rem:GalerkinProjection}.
\par
 The estimate \eqref{eq:EnergyErrors} follows from
 \eqref{eq:Energy1}--\eqref{eq:Energy4} and the triangle inequality.
\end{proof}
\begin{remark}\label{rem:EnergyComparison}\rm
  The estimate \eqref{eq:EnergyErrors} indicates that, up to an
  $O(H+\rho)$ error, the
  approximation of $(\bar y,\bar p)$ by $(\bMSLy,\bMSLp)$ in the
  energy norm is as good as the
  fine scale approximation by $(\byh,\bph)$.
  By comparing \eqref{eq:SmoothEnergy} with \eqref{eq:EnergyErrors},
  we can also say that
  up to the fine scale error
   the performance of the multiscale finite element method (with respect to
   the energy norm) on coarse
   meshes is similar to
  the performance of standard finite element methods for problems with
  smooth coefficients on
  convex domains.
\end{remark}
%%%%%%%%%%%%%%%%%%%%%%%
\section{Numerical Results}\label{sec:Numerics}
 In this section we report the numerical results of two examples, one with
 highly heterogeneous coefficients
  and one with highly oscillatory coefficients.  The domain is the unit
  square $\O=(0,1)\times(0,1)$ for both
  examples, and we use the $Q_1$ element on  uniform rectangular meshes.
  The regularization parameter $\gamma$ is
  taken to be $1$.
\par
 The objective function in our computations is given by
\begin{equation}\label{eq:newJ}
  \tilde J(y,u)=\frac12\big(\|y\|_\LT^2+\gamma\| u\|_\LT^2\big)-\int_\O yy_ddx
\end{equation}
 that differs from $J(y,u)$ by the constant $\|y_d\|_\LT^2/2$.
\par
 The fine scale solution $(\bar y_h,\bar u_h)$ (where $\cT_\rho=\cT_h$) is computed
 by using the primal-dual
 interior point method
 in the PETSc/TAO library with 20 processors on the SuperMIC supercomputer at
 Louisiana State University.  Each compute node
 is equipped with
 two 2.8GHz 10-Core Ivy Bridge-EP E5-2680 Xeon 64-bit Processors, two
 Intel Xeon Phi 7120P Coprocessors,
 64GB DDR3 1866MHz Ram, 500GB HD, 56 Gigabit/sec Infiniband network interface, and
 1 Gigabit Ethernet network interface.
\par
 The DD-LOD solution $(\bMSLy,\buH)$
 (with $\cT_\rho=\cT_H$) is computed
 by using the quadprog algorithm  in MATLAB on a Lenovo Thinkpad X1 Carbon
 laptop
 with a 12th Gen Intel(R) Core(TM) i7-1260P processor,
 4.70 GHz Max Turbo Frequency, an 18MB Intel(R) Smart Cache and
 32 GB of RAM.
\par\vspace{.1in}
\begin{example}[Highly Heterogeneous Coefficients]\label{example:HH}\rm
  The coefficient matrix for this example is given by
\begin{equation*}
  \cA=\begin{bmatrix}
    \cA_{11}& {\bf 0}\\
    {\bf 0} &\cA_{22}
  \end{bmatrix},
\end{equation*}
 where $\cA_{11}$ and $\cA_{22}$ are piecewise constant matrices with respect to a
 $40\times 40$ uniform rectangular subdivision of $\O$.  The values of $\cA_{11}$
 and $\cA_{22}$ on each
  square of the subdivision are randomly generated and range between 1 and 1350
  (cf. Figure~\ref{fig:Heterogeneous}).
\begin{figure}[hhh!]
\begin{center}
  \includegraphics[width=.35\linewidth]{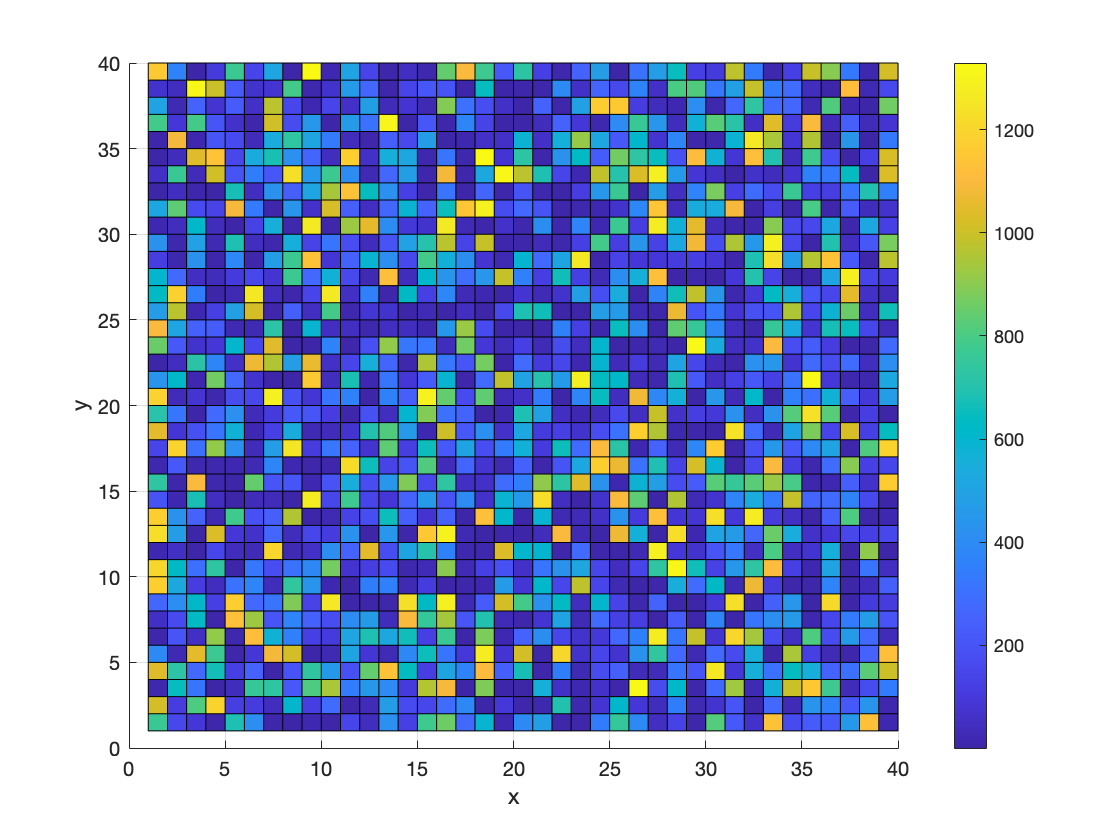} \hspace{20pt}
  \includegraphics[width=.35\linewidth]{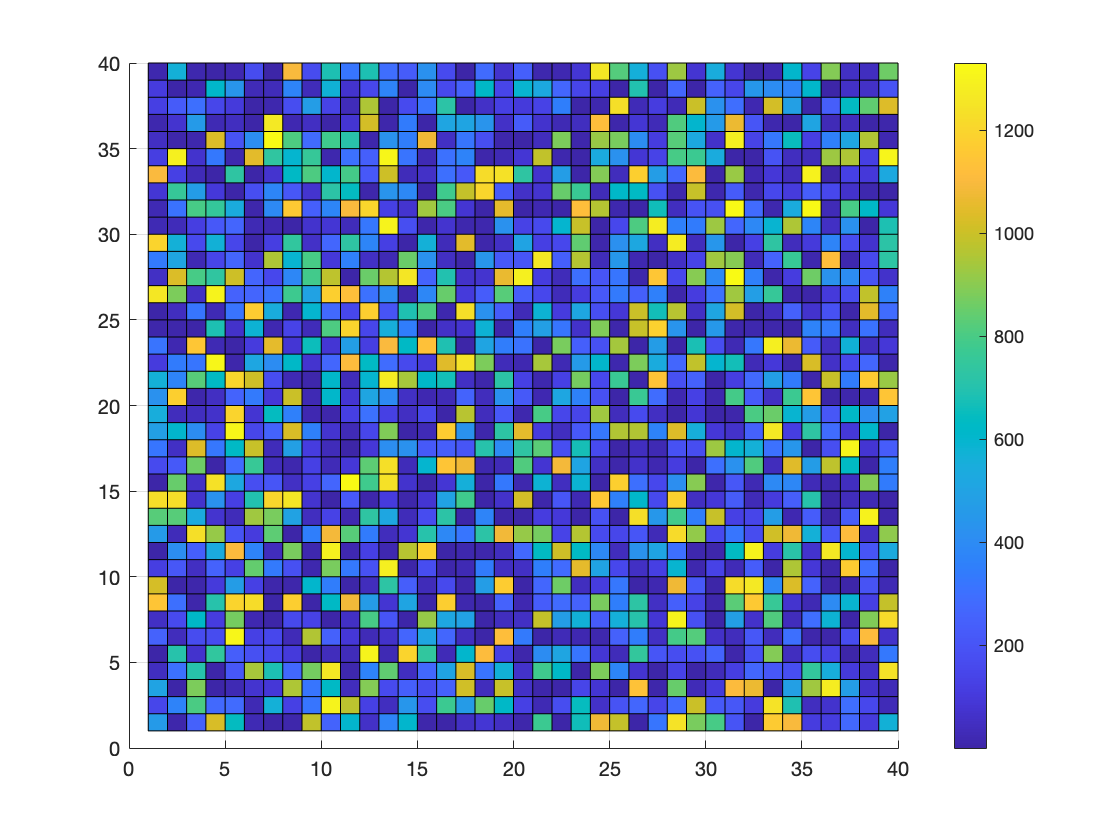}
\end{center}
\caption{$\cA_{11}$ $($left$)$ and $\cA_{22}$ $($right$)$}
\label{fig:Heterogeneous}
\end{figure}
\par
 We choose $y_d=1$ and
 the control constraints are given by
  $\phi_1(x)=0.0002x_1-0.0001$ and $\phi_2(x)=0.0002x_2+0.0001$
 (cf. Figure~\ref{fig:HHCC}).
\begin{figure}[hhh!]
  \centering
  \includegraphics[width=0.3\linewidth]{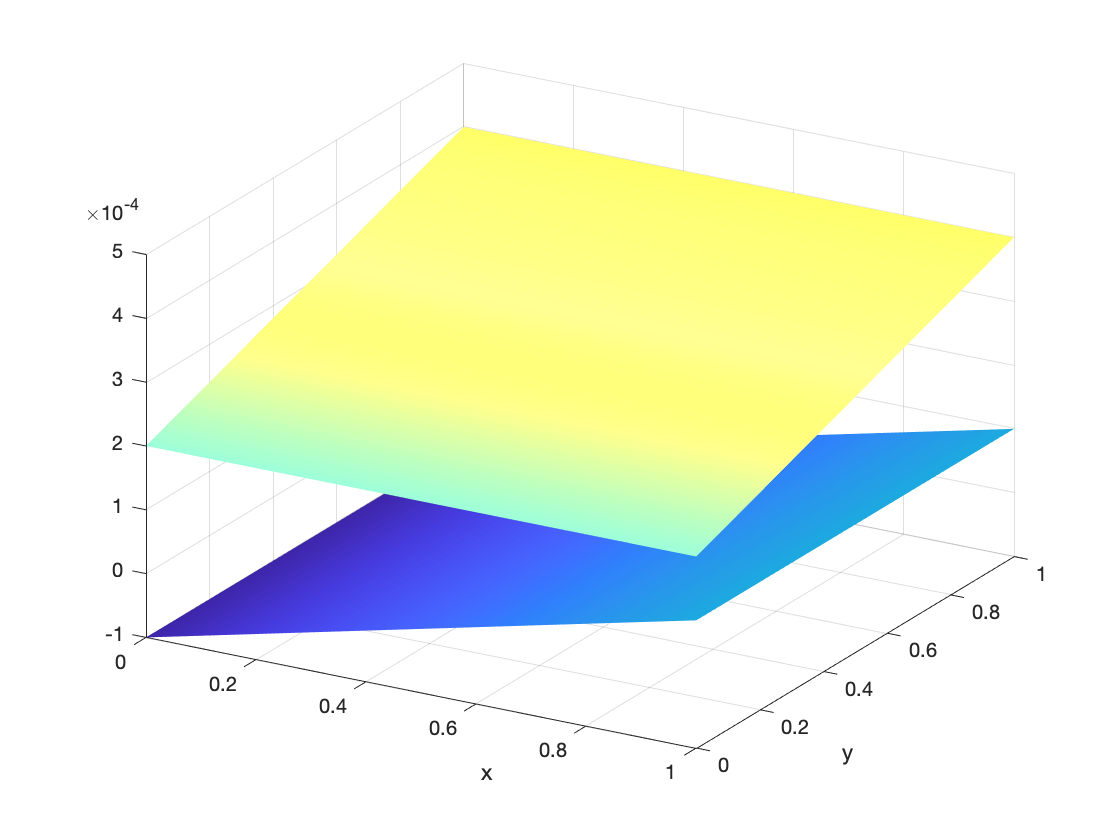}
  \caption{Graphs of the control constraints $\phi_1$ and $\phi_2$ for
  Example~\ref{example:HH}.}
  \label{fig:HHCC}
\end{figure}
\par
 We take $h=1/320$ for the fine scale solution $(\bar y_h,\bar u_h)$.
 In the first set of
 experiments we take
 $H=1/10, 1/20,1/40,1/80$ for the DD-LOD solution $(\bMSLy,\buH)$ with
 $\cT_\rho=\cT_H$.
 The number of iterations $k$ used in the
 solution of the corrector equation equals
 $\lceil-3\ln H\rceil$ for $H=1/10$, $1/20$ and $1/40$, and
 equals $\lceil-6\ln H\rceil$ for $H=1/80$.  The relative errors
 for the approximation of the standard finite element
 solution $(\bar y_h,\bar u_h)$ by the multiscale finite element solution
 $(\bMSLy,\buH)$ are presented in
 Figure~\ref{fig:HHControlStateErrors}.
\begin{figure}[hhh!]
  \centering
\begin{minipage}{2.65in}
\centering
  \includegraphics[width=1\linewidth]{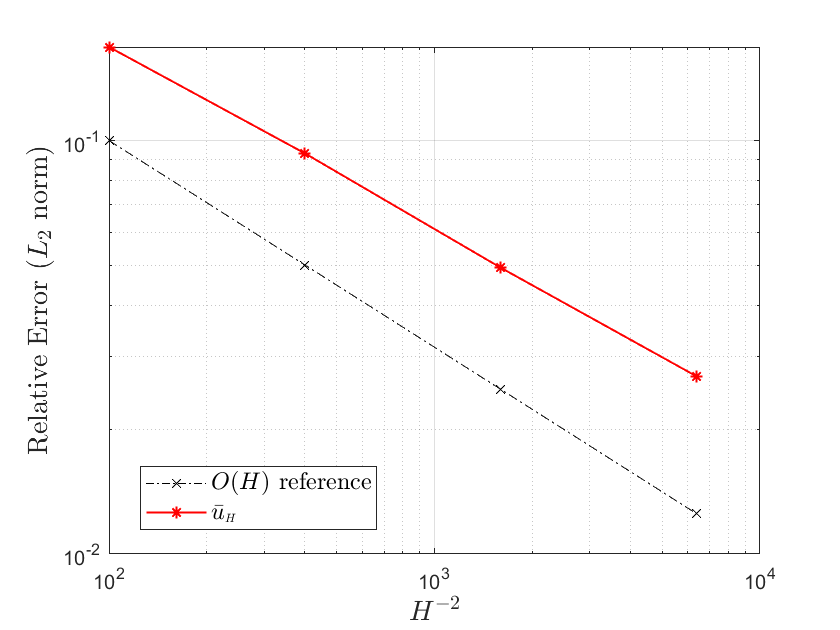}
\par\centerline{\footnotesize(a)}
\end{minipage}\hspace{30pt}
\begin{minipage}{2.65in}
\centering
  \includegraphics[width=1\linewidth]{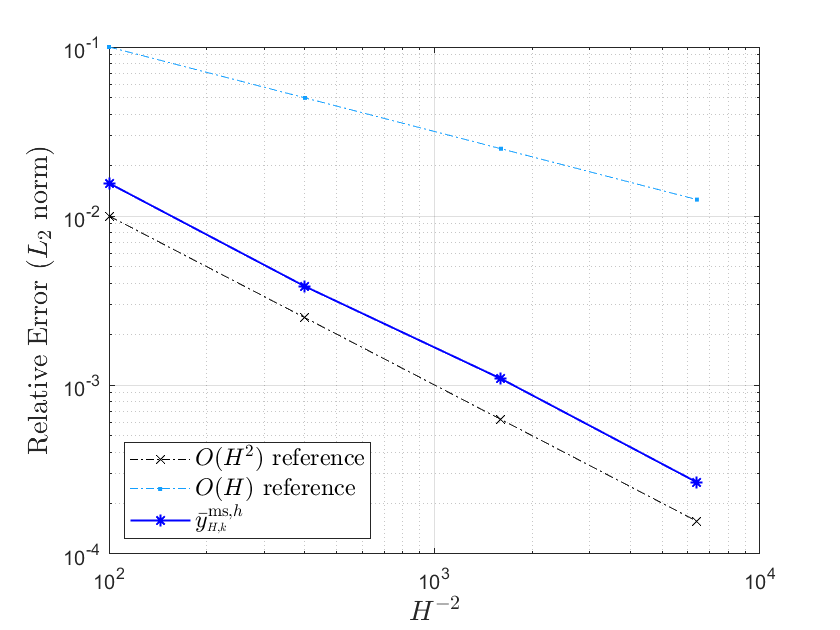}
\par\centerline{\footnotesize(b)}
\end{minipage}
\par
\begin{minipage}{2.65in}
 \centering
  \includegraphics[width=1\linewidth]{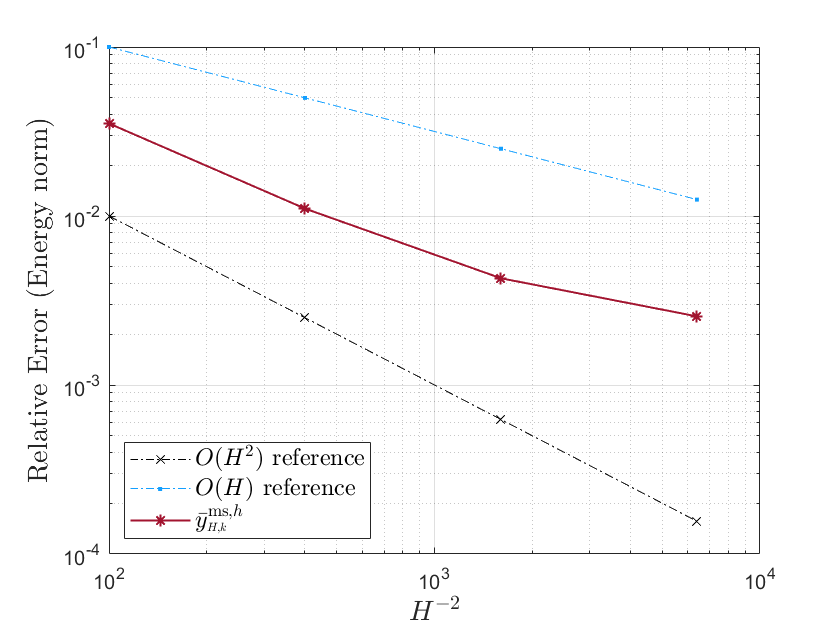}
\par
\centerline{\footnotesize(c)}
\end{minipage}
  \caption{(a) relative $L_2$ error of $\buH$, (b) relative $L_2$ error of  $\bMSLy$
  and (c) relative energy error of $\bMSLy$ for Example~\ref{example:HH} with
  $H=1/10,1/20,1/40,1/80$.}
  \label{fig:HHControlStateErrors}
\end{figure}
\par
 The $O(H)$ convergence of $\buH$ predicted by Theorem~\ref{thm:DDLODError}
 is observed.
 The convergence of $\bMSLy$ in the $L_2$ norm is $O(H^2)$, which is better
 than the $O(H)$ convergence predicted by Theorem~\ref{thm:DDLODError}.
 It should be noted that the error estimate in
 \eqref{eq:DDLODError} concerns the approximation of $(\bar y,\bar u)$ by
  $(\bMSLy,\buH)$, and  the results
 reported in Figure~\ref{fig:HHControlStateErrors} measure the approximation
 of $(\bar y_h,\bar u_h)$ by
 $(\bMSLy,\buH)$.  The convergence of $\bMSLy$ in the energy norm is $O(H)$,
 which agrees with
 Theorem~\ref{thm:EnergyErrors}.
\par
 For this example, the value of the modified cost function $\tilde J$ in
 \eqref{eq:newJ} is $-3.60479\times 10^{-8}$ for the fine scale
 standard finite element
 solution $(\bar y_h,\bar u_h)$.
 The values of $\tilde J(\bMSLy,\buH)$ are displayed in Table~\ref{table:HHMinimum}.
 The order of convergence of $\tilde J(\bMSLy,\buH)$  is roughly $O(H^2)$, which is
 consistent with Theorem~\ref{thm:DDLODError}.
\begin{table}[htb!]
\centering
\begin{tabular}{|c | c|}
\hline&\\[-12pt]
$H$& $\tilde J(\bMSLy,\buH)$\\
\hline &\\[-12pt]
$1/10$&$-3.55321\times10^{-8}$\\
$1/20$&$-3.59107\times10^{-8}$\\
$1/40$&$-3.60102\times 10^{-8}$\\
$1/80$&$-3.60372\times 10^{-8}$\\
 \hline
\end{tabular}
\par\medskip
\caption{Values of $\tilde J(\bMSLy,\buH)$ for Example~\ref{example:HH}.}
\label{table:HHMinimum}
\end{table}
 \par
  We compare the graphs of $\bar y_h$ and $\bMSLy$ (with $H=1/20$)
  in Figure~\ref{fig:HHStateComparison},
 and the graphs of $\bar u_h$ and $\buH$ (with $H=1/20$) in
 Figure~\ref{fig:HHControlComparison}.
\begin{figure}[hhh!]
  \centering
 \includegraphics[width=0.3\linewidth]{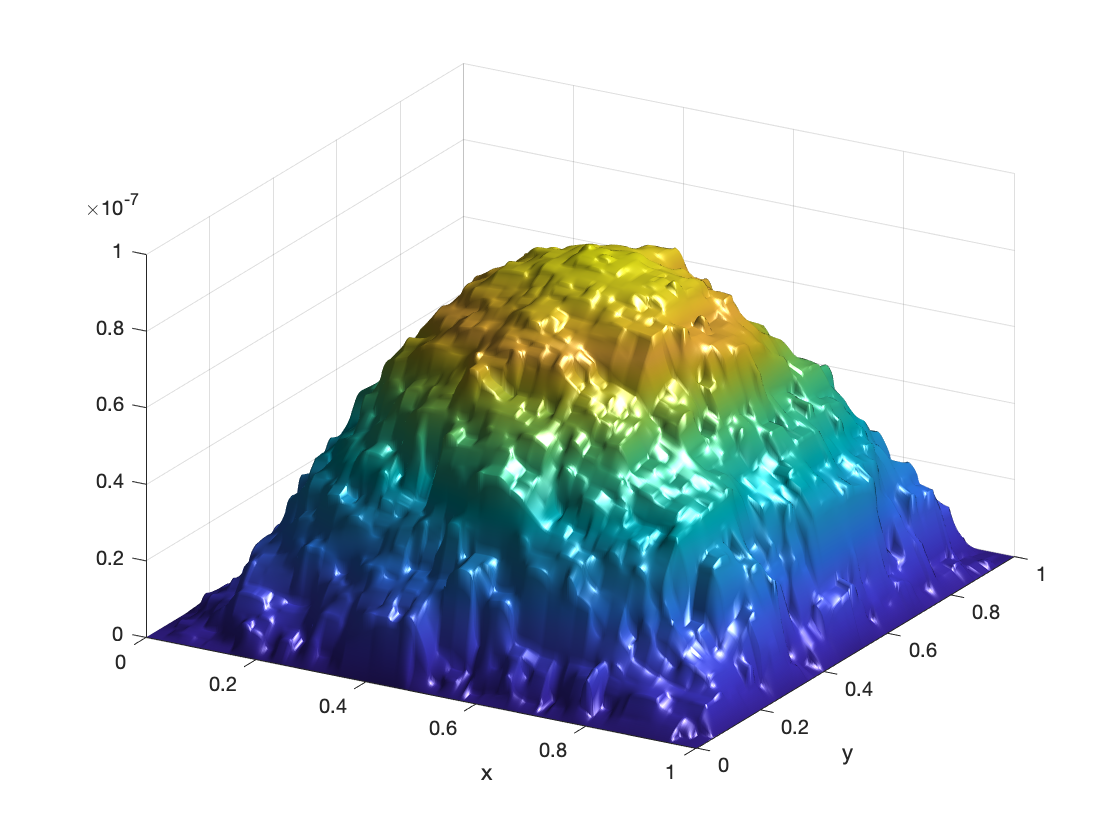}
  \hspace{20pt}
  \includegraphics[width=0.3\linewidth]{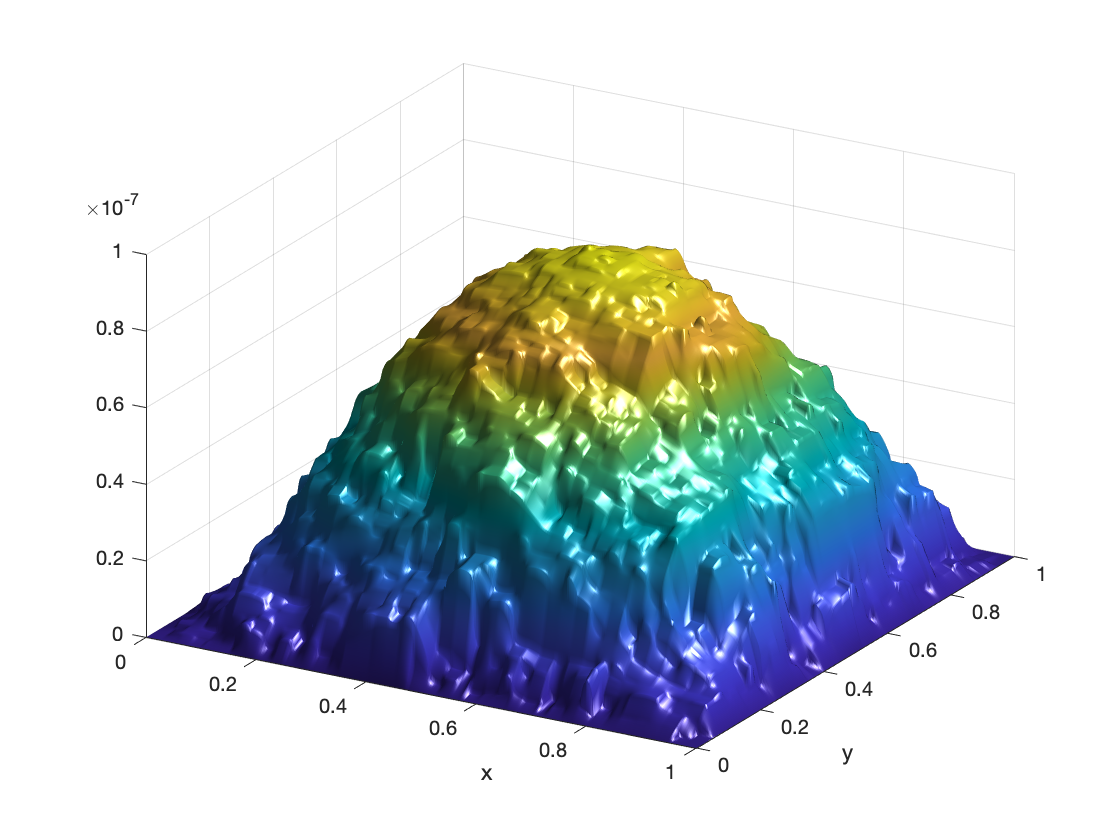}
\caption{Graph of $\bar y_h$ (left) and graph of $\bMSLy$ (right, with $H=1/20$)
for Example~\ref{example:HH}.}
\label{fig:HHStateComparison}
\end{figure}
\begin{figure}[hhh!]
  \centering
 \includegraphics[width=0.3\linewidth]{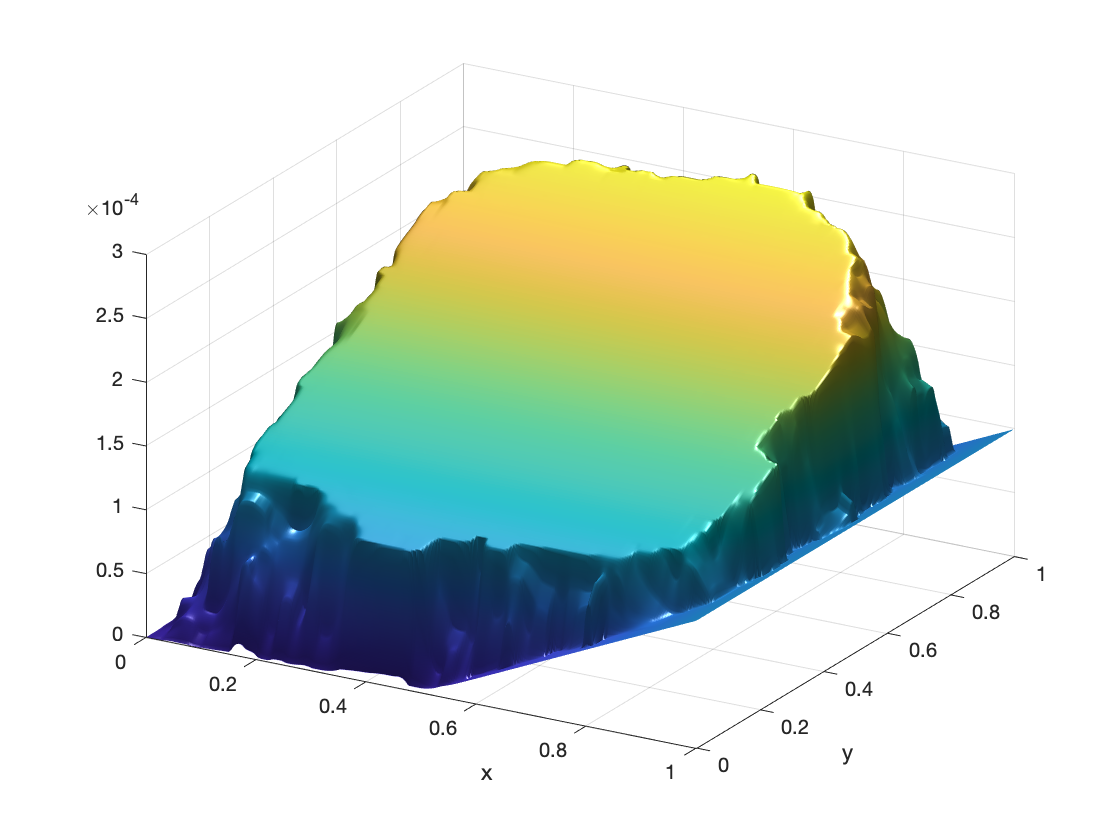}
\hspace{20pt}
 \includegraphics[width=0.3\linewidth]{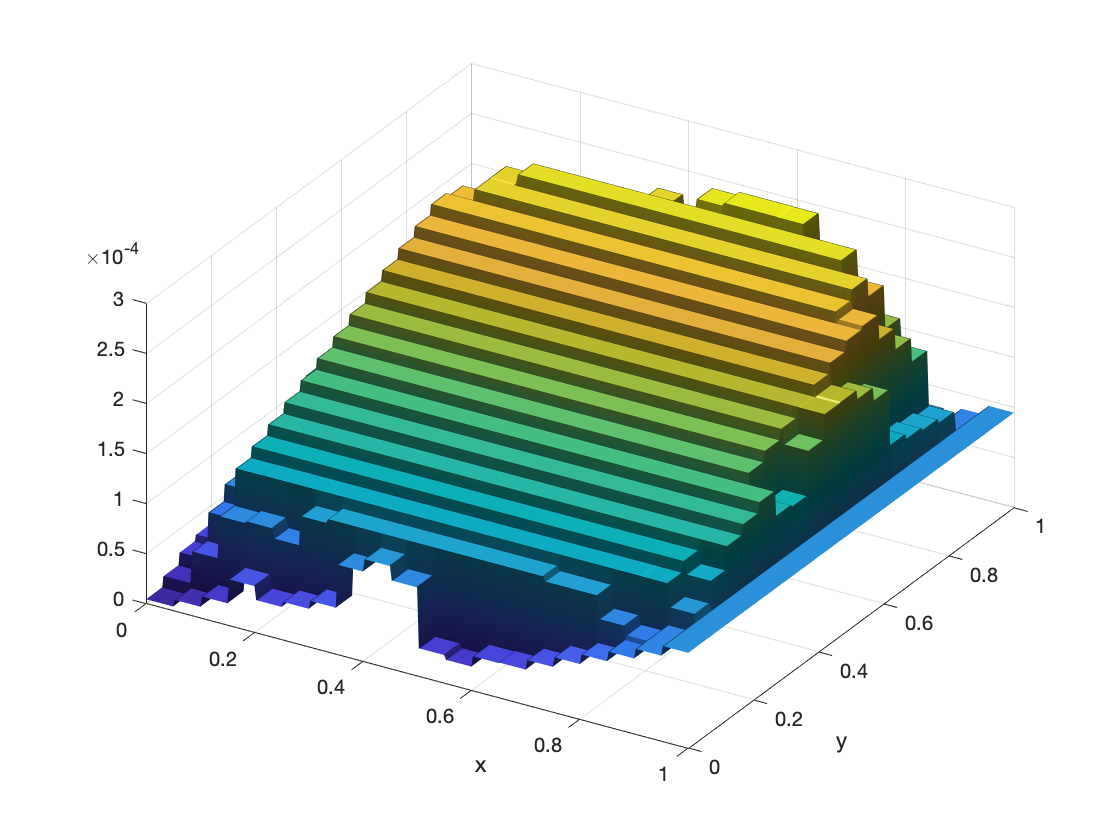}
\caption{Graph of $\bar u_h$ (left) and graph of $\buH$ (right, with $H=1/20$)
for Example~\ref{example:HH}.}
\label{fig:HHControlComparison}
\end{figure}
\par\phantom{a}
\par\vspace{-10pt}\indent
 The active sets for $\bar u_h$ and $\buH$ (with $H=1/20$) are depicted in
 Figure~\ref{fig:HHActiveSet1}  and Figure~\ref{fig:HHActiveSet2}.
\par
 \begin{figure}[hhh!]
  \centering
  \includegraphics[width=0.3\linewidth]{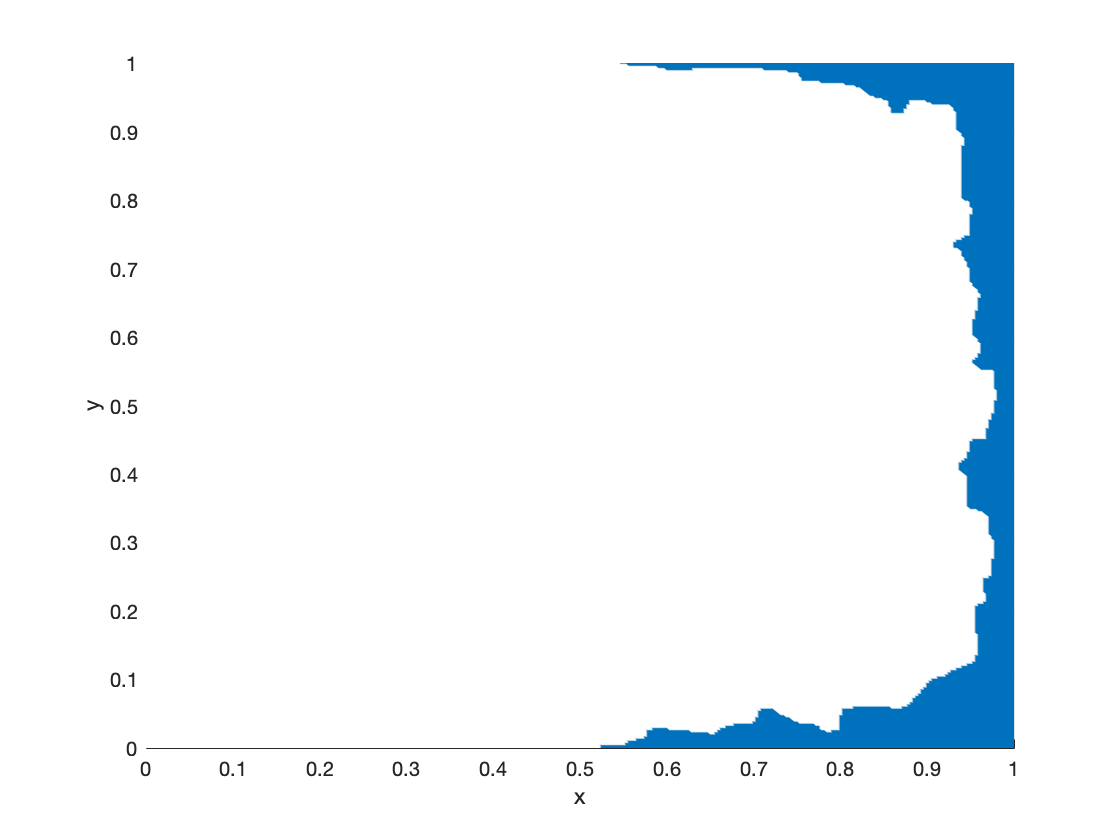}
\hspace{20pt}
  \includegraphics[width=0.3\linewidth]{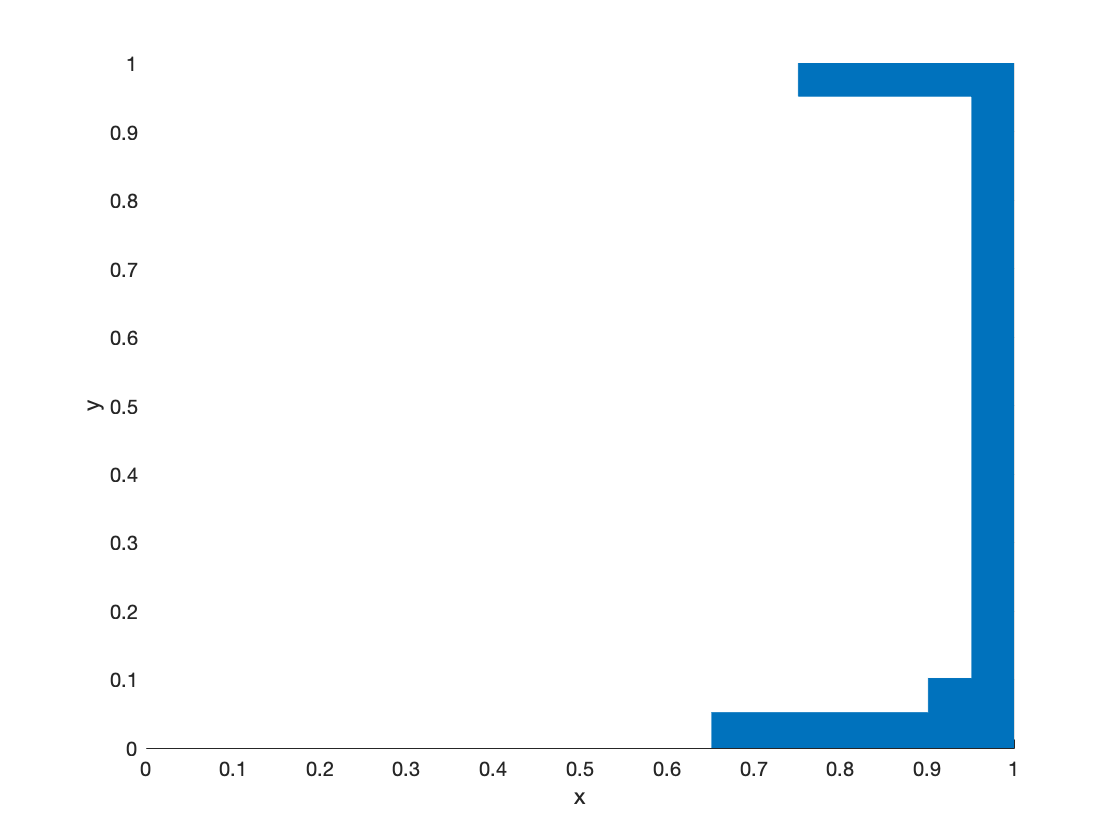}
\caption{Active sets for $\phi_1$ for Example~\ref{example:HH}:
 $\bar u_h$ (left) and
 $\buH$ (right, $H=1/20$). }
\label{fig:HHActiveSet1}
\end{figure}
\begin{figure}[hhh!]
\centering
  \includegraphics[width=0.3\linewidth]{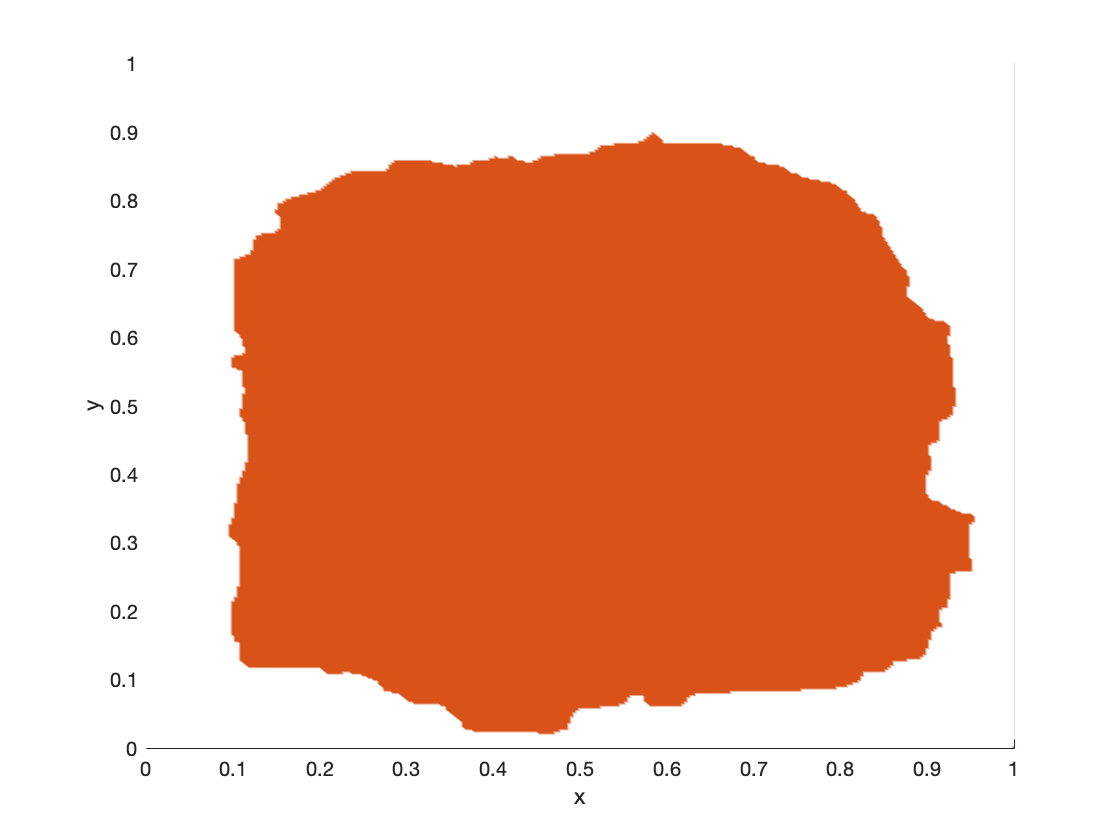}
\hspace{20pt}
  \includegraphics[width=0.3\linewidth]{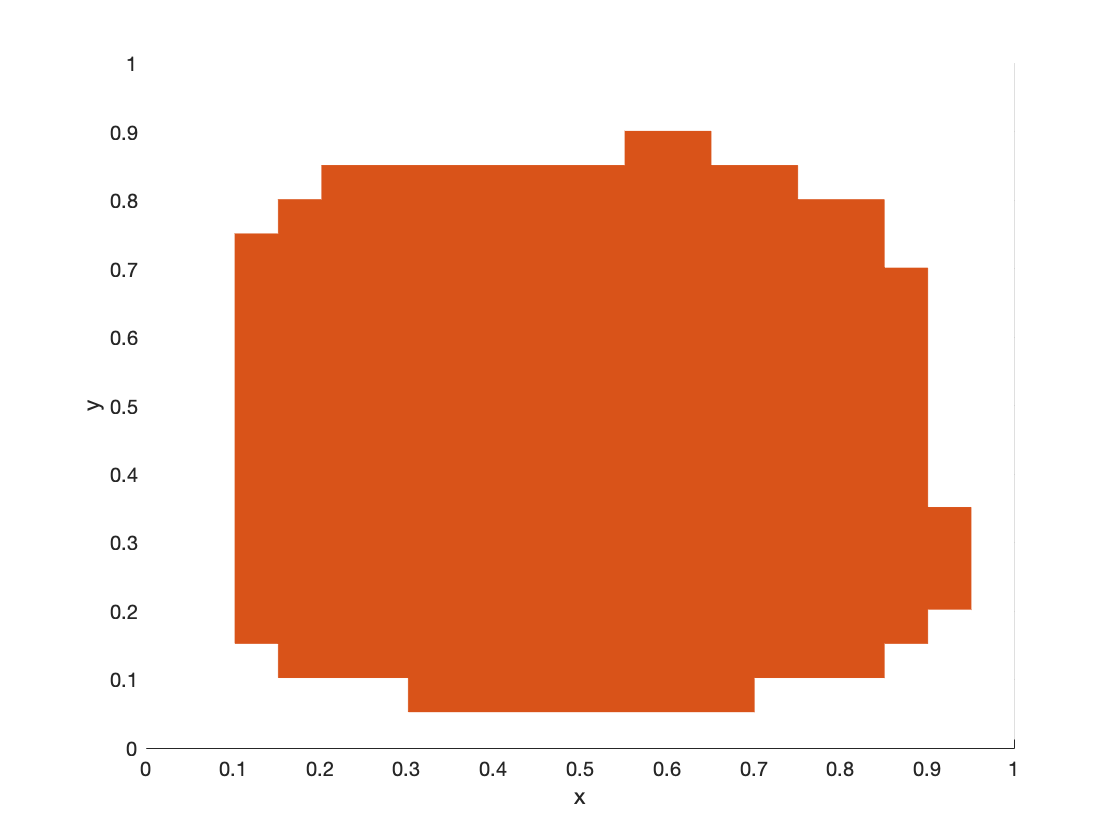}
\caption{Active sets for $\phi_2$ for Example~\ref{example:HH}:
$\bar u_h$ (left) and
 $\buH$ (right, $H=1/20$).}
\label{fig:HHActiveSet2}
\end{figure}
\par
 The computation of the fine scale standard finite element solution
 $(\bar y_h,\bar u_h)$ of the discrete optimization problem takes
 $3.90\times 10^{+1}$
 seconds by using the PETSc/TAO
 library with 20 processors.
 The computational time (in seconds) for $(\bMSLy,\buH)$  using MATLAB
 on a laptop
 are presented in Table~\ref{table:HHTimes} for $H=1/10,1/20,1/40$.
\begin{table}[htb!]
\centering
\begin{tabular}{|c|c|}
  \hline
  $H$ & Time\\
  \hline&\\[-12pt]
  $1/10$ & $1.26\times 10^{-2}$\\
  $1/20$ & $1.74\times 10^{-1}$\\
  $1/40$ & $1.04\times 10^{+1}$\\
  \hline
\end{tabular}
\par\medskip
\caption{Computational time in seconds for $(\bMSLy,\buH)$
 (Example~\ref{example:HH}).}
\label{table:HHTimes}
\end{table}
\par
 For $H=1/20$, the DD-LOD solution $(\bMSLy,\buH)$ yields a reasonable
 approximation of $(\bar y_h,\bar u_h)$
 (cf. Figures~\ref{fig:HHStateComparison}--\ref{fig:HHActiveSet2}) and
 its computation is more than
  100 times faster
 than the computation of $(\bar y_h,\bar u_h)$.
\par
 In the second set of  experiments we take $H=1/20$ and
 $\rho=1/40, 1/80, 1/160$ for the
 DD-LOD solution $(\bMSLy,\bar u_{\rho})$.  In view of
 Theorem~\ref{thm:DDLODError}
 and Theorem~\ref{thm:EnergyErrors},
 we expect these approximate solutions will improve over the
 approximate solution
 $(\bMSLy,\buH)$ with $H=1/20$
 and $\cT_\rho=\cT_H$
 obtained in the first set of experiments.  This is confirmed by
 comparing the values of the cost function
 $\tilde J$ in Table~\ref{table:HHrhoMinimum} with the value
  $\tilde J(\bar y_h,\bar u_h)=-3.60479\times 10^{-8}$ for
 the fine scale solution.  The number of significant digits increases
 from $2$ to $4$ as
 $\rho$ decreases from $1/20$ to $1/160$.
\begin{table}[htb!]
\centering
\begin{tabular}{|c | c|}
\hline&\\[-12pt]
$\rho$& $\tilde J(\bMSLy,\bar u_{\rho})$\\
\hline &\\[-12pt]
$1/20$&$-3.59107\times10^{-8}$\\
$1/40$&$-3.60090\times 10^{-8}$\\
$1/80$&$-3.60357\times 10^{-8}$\\
$1/160$&$-3.60431\times 10^{-8}$\\
 \hline
\end{tabular}
\par\medskip
\caption{Values of $\tilde J(\bMSLy,\bar u_\rho)$ with
 $H=1/20$ and various $\rho$  for Example~\ref{example:HH}.}
\label{table:HHrhoMinimum}
\end{table}
\par
 We can also visualize the improvement due to a smaller $\rho$ by comparing
 the graph of the
 fine scale solution $\bar u_h$ for the optimal control and the graph of
 the DD-LOD solution
 $\bar u_\rho$ for the optimal control
 (with $H=1/20$ and $\rho=1/160$) in Figure~\ref{fig:HHrhoControlComparison}.
 They are hardly distinguishable, which is not the case for the  graphs in
  Figure~\ref{fig:HHControlComparison}.
\begin{figure}[htb!]
  \centering
 \includegraphics[width=0.3\linewidth]{OptControl-STD-319-hh-AL}
\hspace{20pt}
 \includegraphics[width=0.3\linewidth]{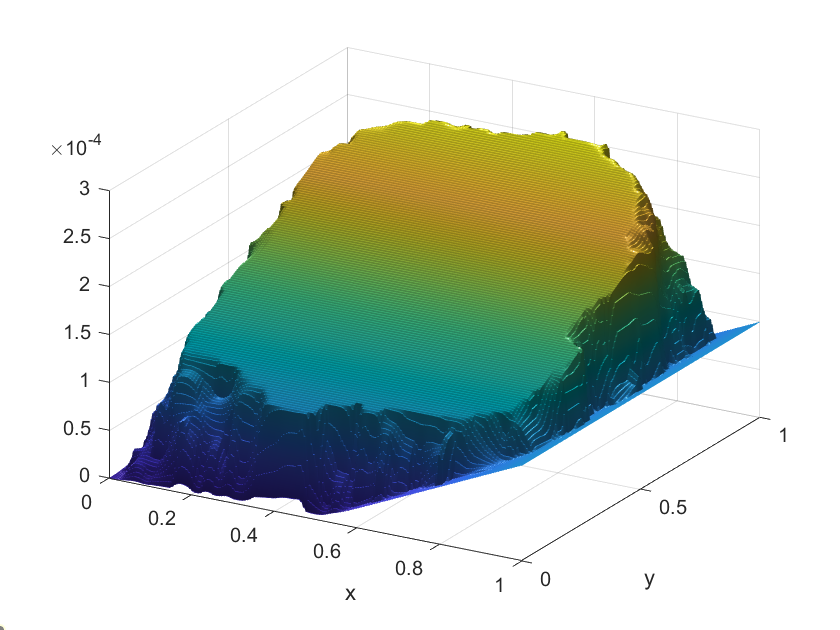}
\caption{Graph of $\bar u_h$ (left) and graph of $\bar u_\rho$
 (right, with $H=1/20$ and $\rho=1/160$.)
for Example~\ref{example:HH}.}
\label{fig:HHrhoControlComparison}
\end{figure}
\par
 This is also true for the active sets, where the ones for the fine scale
 solution $\bar u_h$
 and the ones for the DD-LOD solution  $\bar u_\rho$ (with $H=1/20$ and
 $\rho=1/160$)
  are almost identical in
 Figure~\ref{fig:HHrhoActiveSet1} and Figure~\ref{fig:HHrhoActiveSet2}.
 \begin{figure}[hhh!]
  \centering
  \includegraphics[width=0.3\linewidth]{Active-set-part1-af-19-319-hh}
\hspace{20pt}
  \includegraphics[width=0.3\linewidth]{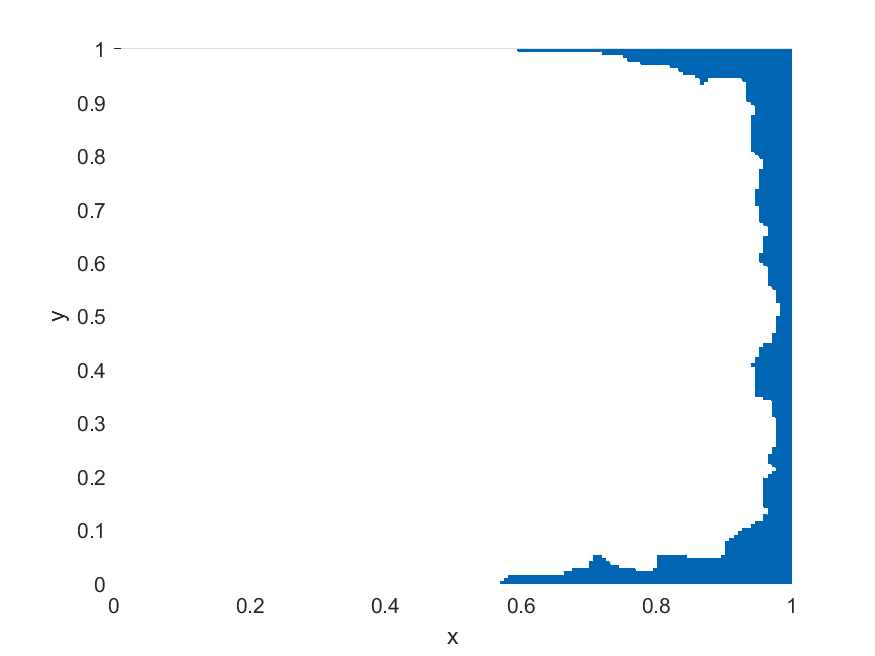}
\caption{Active sets for $\phi_1$ for Example~\ref{example:HH}:
 $\bar u_h$ (left) and
 $\bar u_\rho$ (right, $H=1/20$ and $\rho=1/160$). }
\label{fig:HHrhoActiveSet1}
\end{figure}
\begin{figure}[hhh!]
\centering
  \includegraphics[width=0.3\linewidth]{Active-set-part2-af-19-319-hh}
\hspace{20pt}
  \includegraphics[width=0.3\linewidth]{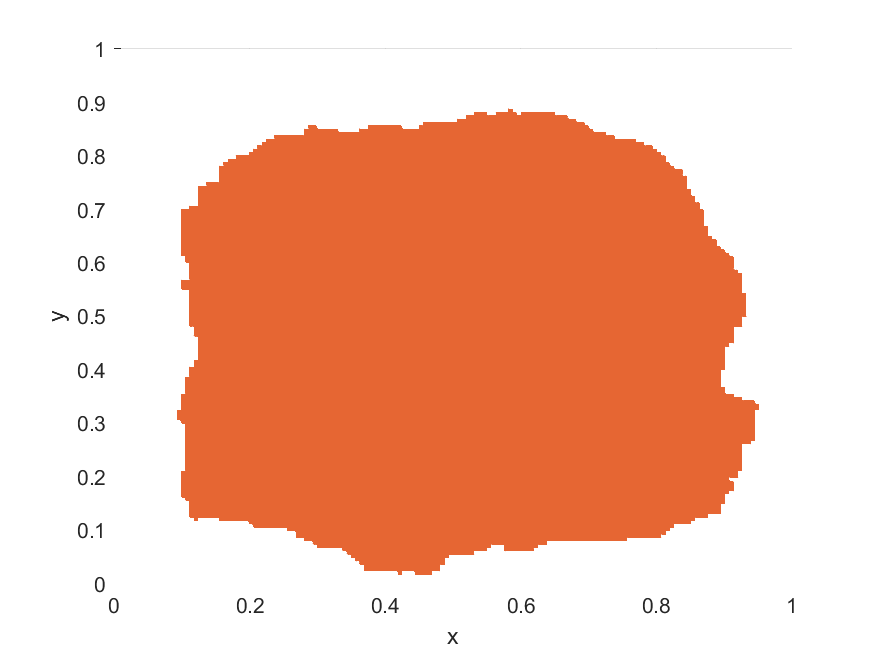}
\caption{Active sets for $\phi_2$ for Example~\ref{example:HH}:
 $\bar u_h$ (left) and
 $\bar u_\rho$ (right, $H=1/20$ and $\rho=1/160$).}
\label{fig:HHrhoActiveSet2}
\end{figure}
\end{example}
%
%%%%%%%%%%%%%%%%%%%%%%%%%%%%%%
\begin{example}[Highly Oscillatory Coefficients]\label{example:HO}\rm
    The coefficient matrix for this example is given by
\begin{equation*}
  \cA=\begin{bmatrix}
    c(x)& {\bf 0}\\
    {\bf 0} & c(x)
  \end{bmatrix},
\end{equation*}
 where
\begin{equation*}
  c(x)=\frac{2+1.8\sin\left(\frac{2\pi x_1}{\epsilon}\right)}
  {2+1.8\sin\left(\frac{2\pi x_2}{\epsilon}\right)}+
  \frac{2+\sin\left(\frac{2\pi x_2}{\epsilon}\right)}
  {2+1.8\sin\left(\frac{2\pi x_1}{\epsilon}\right)}
\end{equation*}
 with $\epsilon=0.025$.  This choice of coefficients originates
 from the pioneering work \cite{HW:1999:MS} in numerical homogenization.
\par
  We choose $y_d=-1$ and the control constraints are given by
$\phi_1(x)=-0.01x_1-0.005$ and $\phi_2(x)=0.0007x_2-0.005$
 (cf. Figure~\ref{fig:HOCC}).
\begin{figure}[hhh!]
  \centering
  \includegraphics[width=0.4\linewidth]{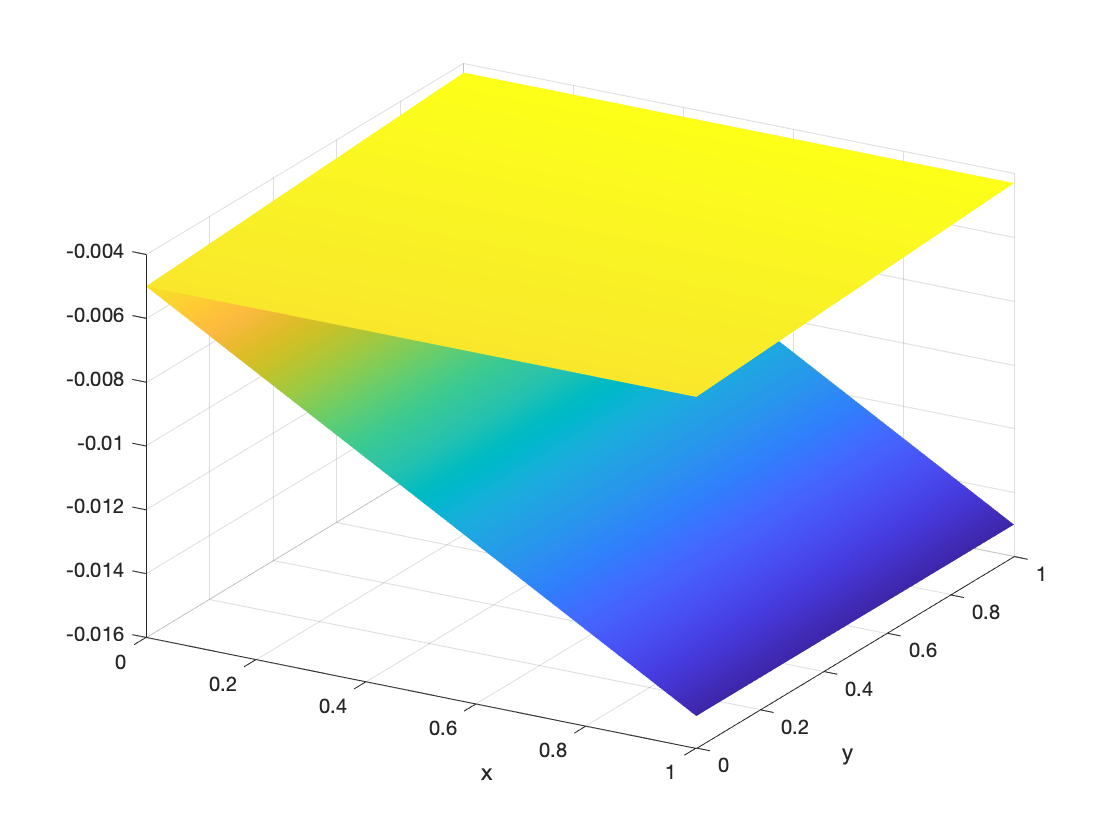}
  \caption{Graphs of the control constraints $\phi_1$ and $\phi_2$ for
   Example~\ref{example:HO}. }
  \label{fig:HOCC}
\end{figure}
\par
  We take $h=1/320$ for the fine scale solution $(\bar y_h,\bar u_h)$.
  In the first set of experiments we compute the DD-LOD solution
  $(\bMSLy,\buH)$ for $H=1/10, 1/20,1/40,1/80$ (with $\cT_\rho=\cT_H$).
  The number of iterations $k$ used in the
 solution of the corrector equation equals  $\lceil-3\ln H\rceil$ for all $H$.
   The relative errors for the approximation of the fine scale standard
   finite element
 solution $(\bar y_h,\bar u_h)$ by the multiscale finite element
 solution $(\bMSLy,\buH)$ are presented in
 Figure~\ref{fig:HOControlStateErrors}.
 The $O(H)$ convergence is observed for both $\buH$ and $\bMSLy$,
 which agrees with Theorem~\ref{thm:DDLODError}
 and Theorem~\ref{thm:EnergyErrors}.
\begin{figure}[hhh!]
  \centering
\begin{minipage}{2.7in}
\centering
  \includegraphics[width=\linewidth]{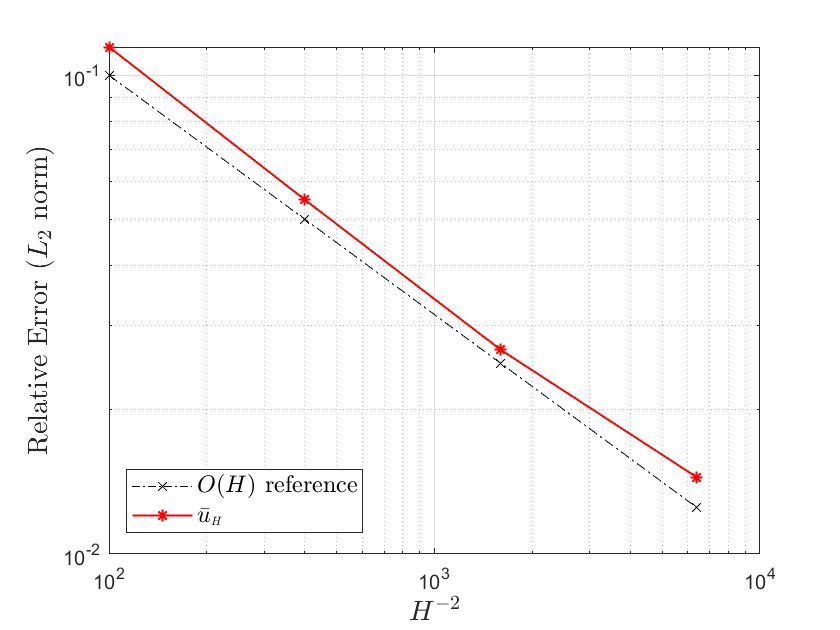}
\par\centerline{\footnotesize (a)}
\end{minipage}\hspace{30pt}
\begin{minipage}{2.7in}
  \includegraphics[width=\linewidth]{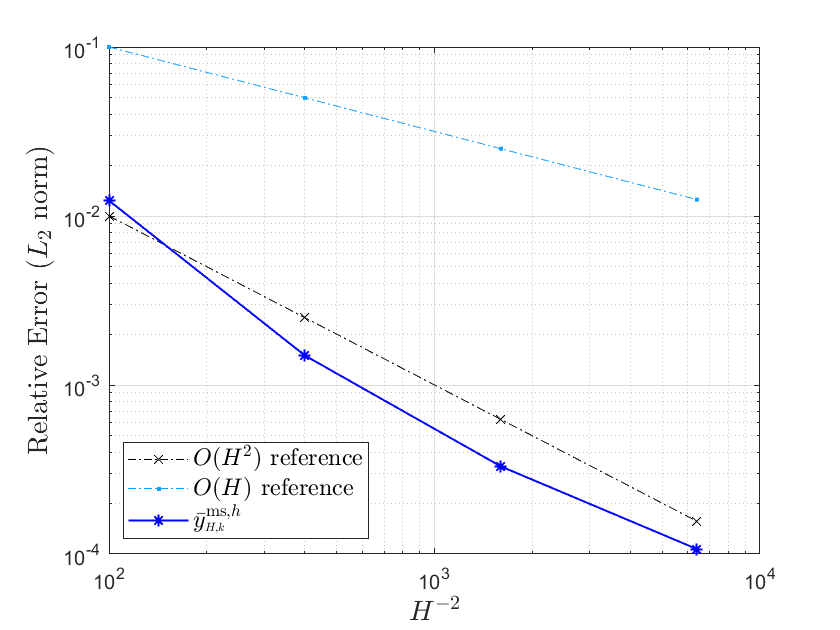}
\par\centerline{\footnotesize (b)}
\end{minipage}
\par
\begin{minipage}{2.7in}
\centering\includegraphics[width=\linewidth]{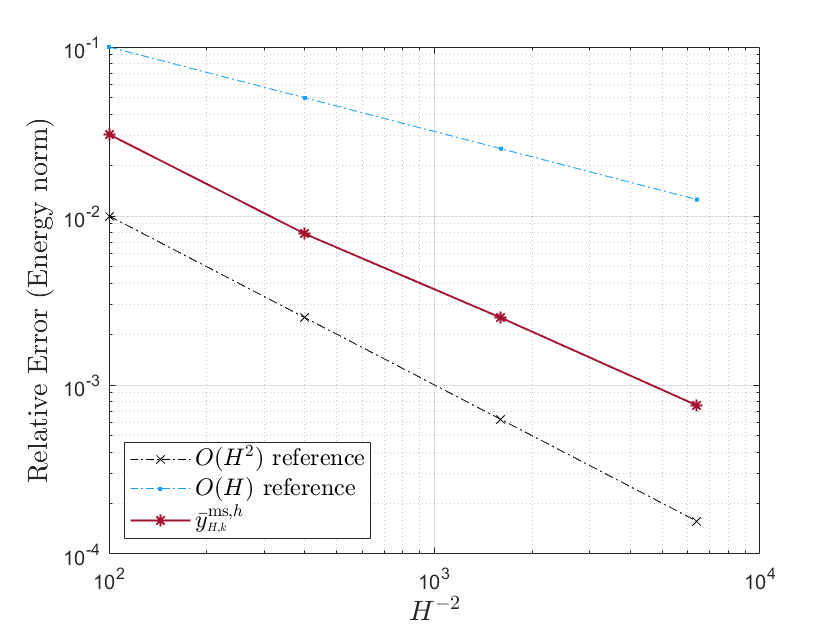}
\par\centerline{\footnotesize (c)}
\end{minipage}
  \caption{(a) relative $L_2$ error of $\buH$, (b) relative $L_2$ error of $\bMSLy$
  and (c) relative energy error of $\bMSLy$ for Example~\ref{example:HO}
   with $H=1/10,1/20,1/20,1/80$.}
  \label{fig:HOControlStateErrors}
\end{figure}
\par
 For this example, the value of the modified cost function $\tilde J$ in
 \eqref{eq:newJ} is $-8.29631\times 10^{-5}$ for the fine scale standard
 finite element solution  $(\bar y_h,\bar u_h)$.
 The values of $\tilde J(\bMSLy,\buH)$ are displayed in Table~\ref{table:HOMinimum}.
 The $O(H^2)$ convergence of $\tilde J(\bMSLy,\buH)$  also agrees with
 Theorem~\ref{thm:DDLODError}.
\begin{table}[htb!]
\centering
\begin{tabular}{|c | c|}
\hline &\\[-12pt]
$H$& $\tilde J(\bMSLy,\buH)$\\
\hline &\\[-12pt]
$1/10$&$-8.22171\times10^{-5}$\\
$1/20$&$-8.28313\times10^{-5}$\\
$1/40$&$-8.29343\times10^{-5}$\\
$1/80$&$-8.29550\times10^{-5}$\\
 \hline
\end{tabular}
\par\medskip
\caption{Values of $\tilde J(\bMSLy,\buH)$ for Example~\ref{example:HO}.}
\label{table:HOMinimum}
\end{table}
\par
 We compare the graphs of $\bar y_h$ and $\bMSLy$ (with $H=1/20$) in
 Figure~\ref{fig:HOStateComparison},
 and the graphs of $\bar u_h$ and $\buH$ (with $H=1/20$)
 in Figure~\ref{fig:HOControlComparison}.
 The active sets for $\bar u_h$ and $\buH$ (with $H=1/20$) are depicted in
 Figure~\ref{fig:HOActiveSet1}
 and Figure~\ref{fig:HOActiveSet2}.
\begin{figure}[hhh!]
\centering
  \includegraphics[width=0.3\linewidth]{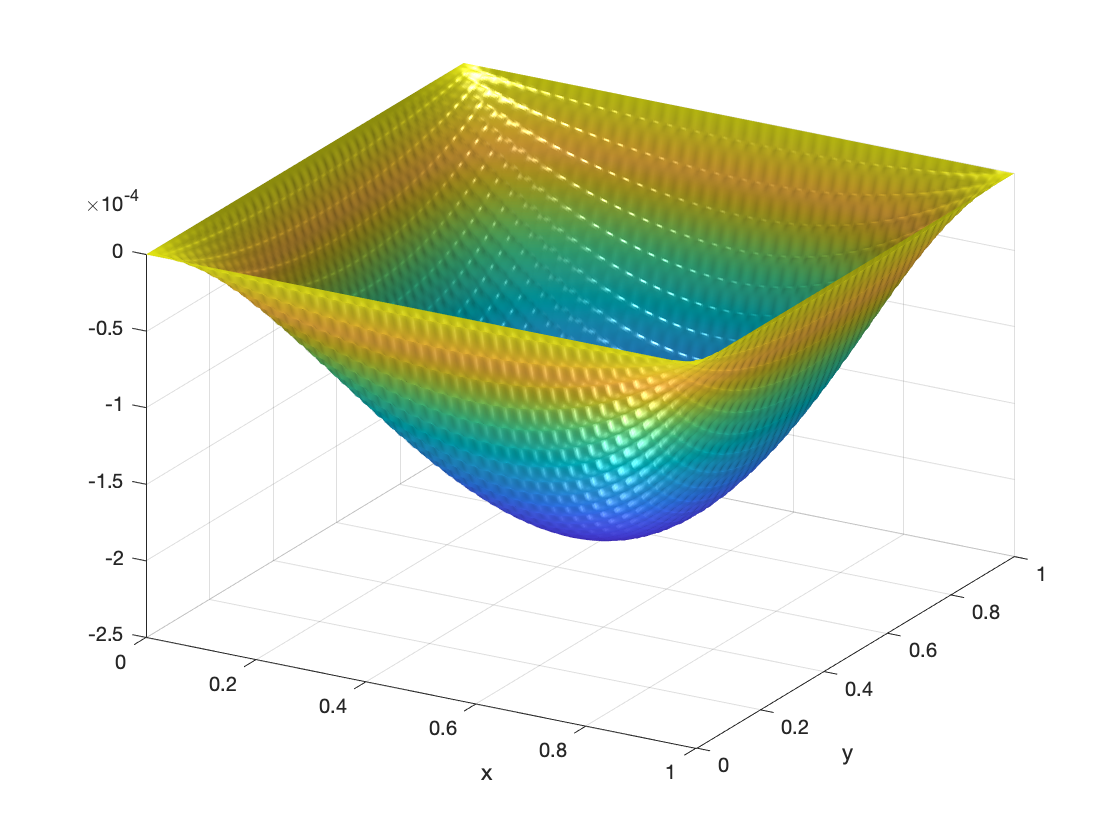}
\hspace{20pt}
  \includegraphics[width=0.3\linewidth]{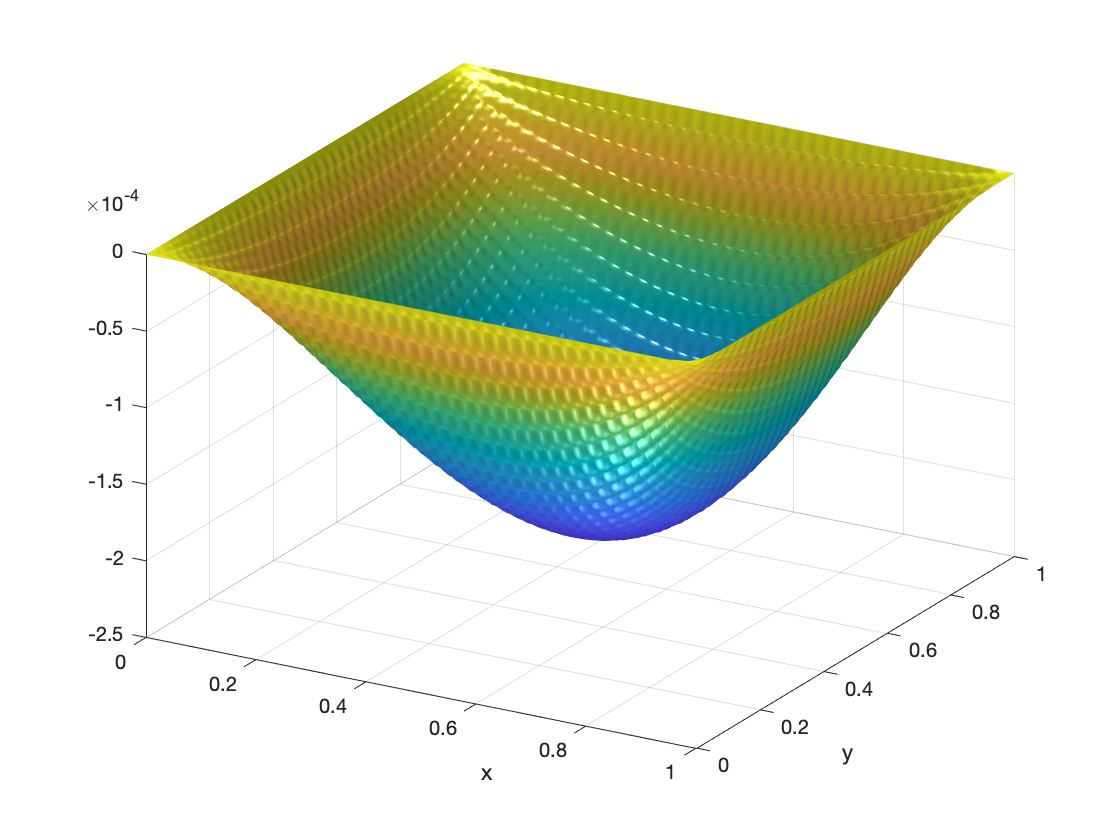}
\caption{Graph of $\bar y_h$ (left) and graph of $\bMSLy$ (right, $H=1/20$)
for Example~\ref{example:HO}.}
\label{fig:HOStateComparison}
\end{figure}
 \begin{figure}[hhh!]
\centering
  \includegraphics[width=0.3\linewidth]{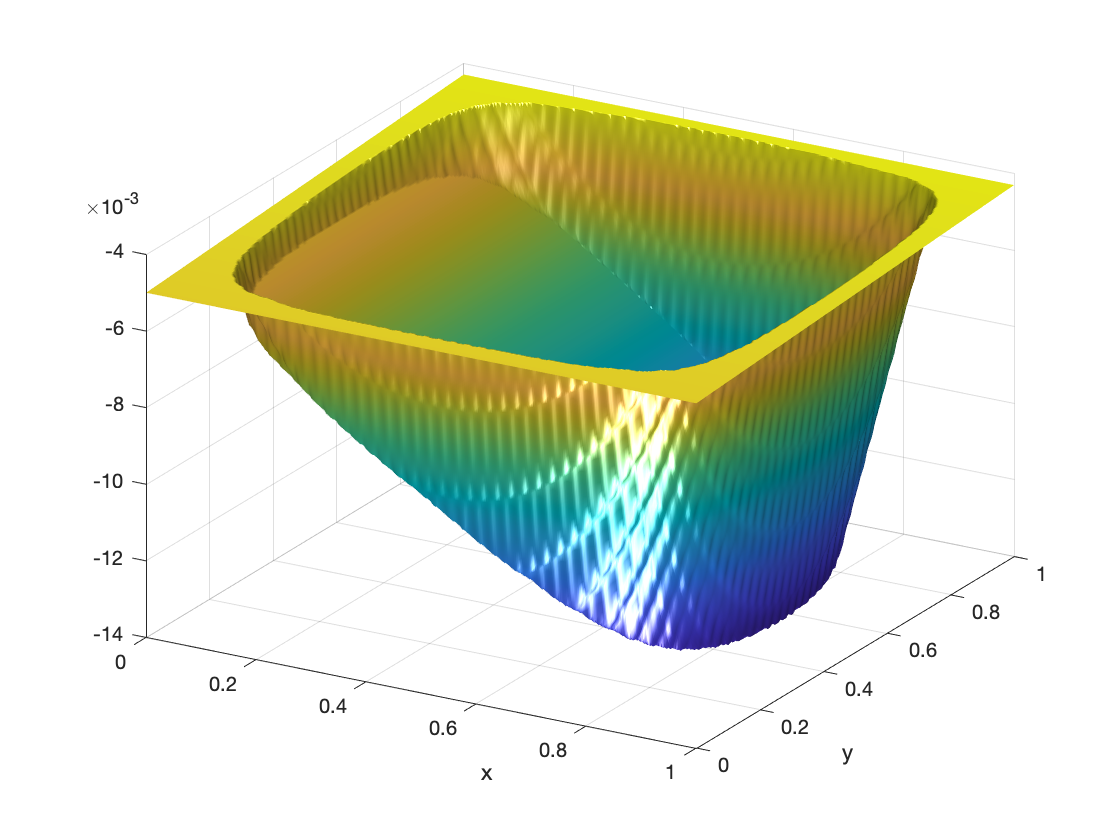}
\hspace{20pt}
  \includegraphics[width=0.3\linewidth]{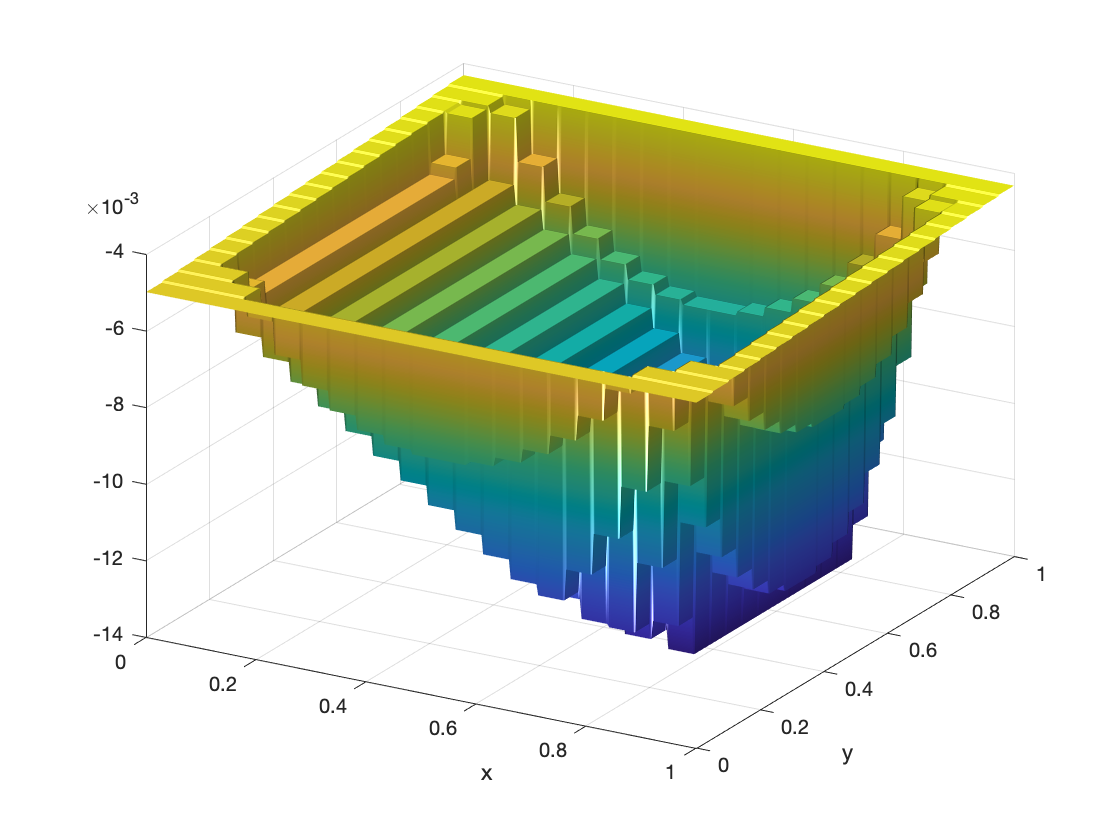}
\caption{Graph of $\bar u_h$ (left) and graph of $\buH$ (right, $H=1/20$)
for Example~\ref{example:HO}.}
\label{fig:HOControlComparison}
\end{figure}
\begin{figure}[hhh!]
\centering
  \includegraphics[width=0.3\linewidth]{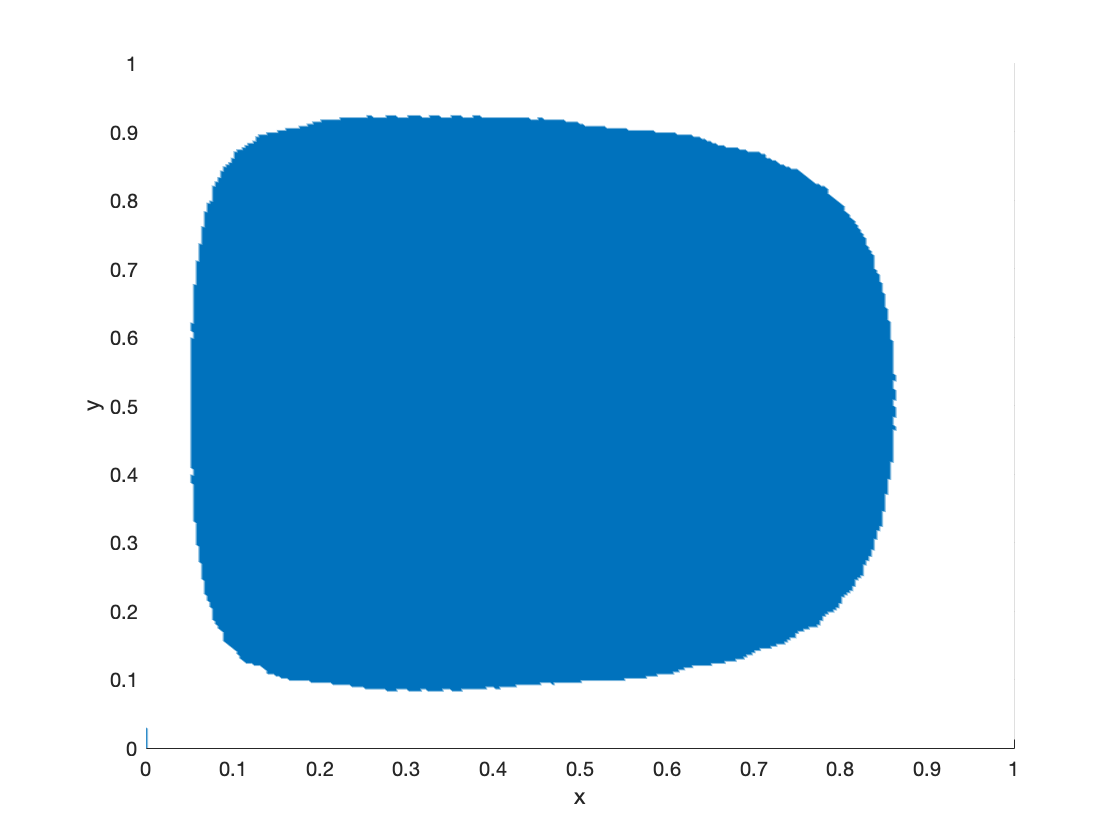}
\hspace{20pt}
  \includegraphics[width=0.3\linewidth]{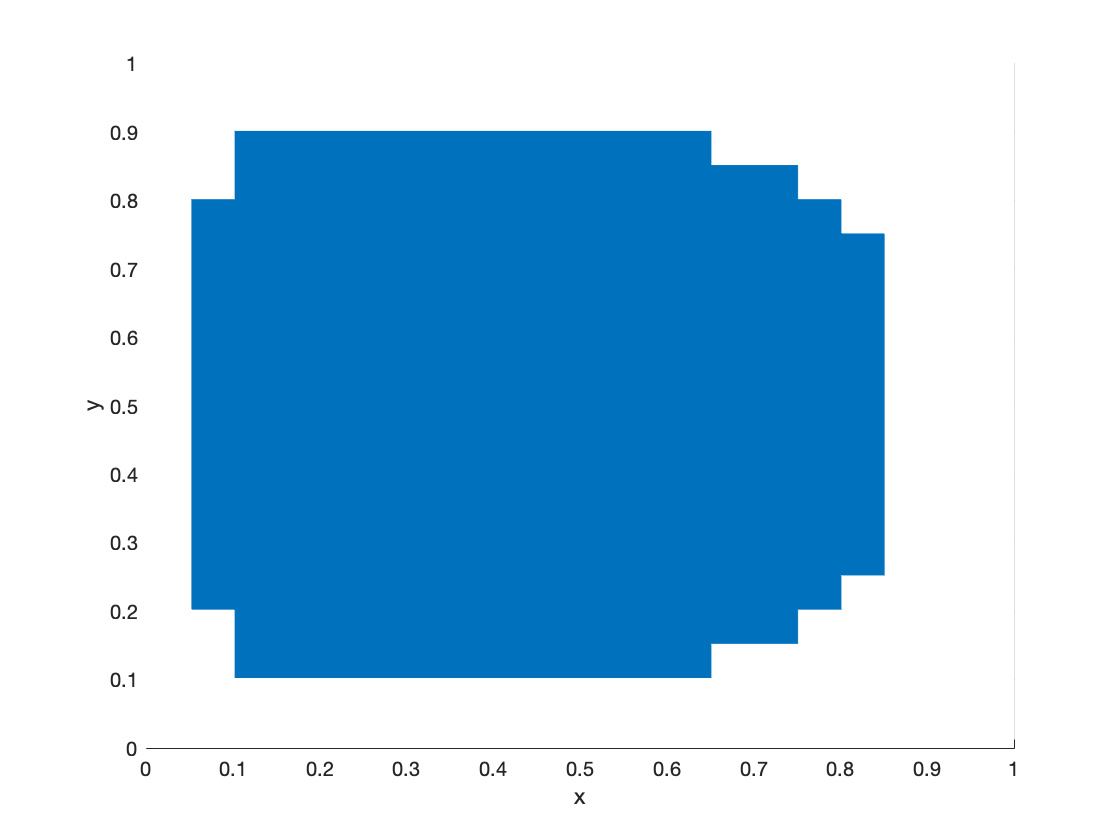}
\caption{Active set for $\phi_1$ for Example~\ref{example:HO}: $\bar u_h$ (left) and
$\buH$ (right, $H=1/20$).}
\label{fig:HOActiveSet1}
\end{figure}
\begin{figure}[hhh!]
\centering
  \includegraphics[width=0.3\linewidth]{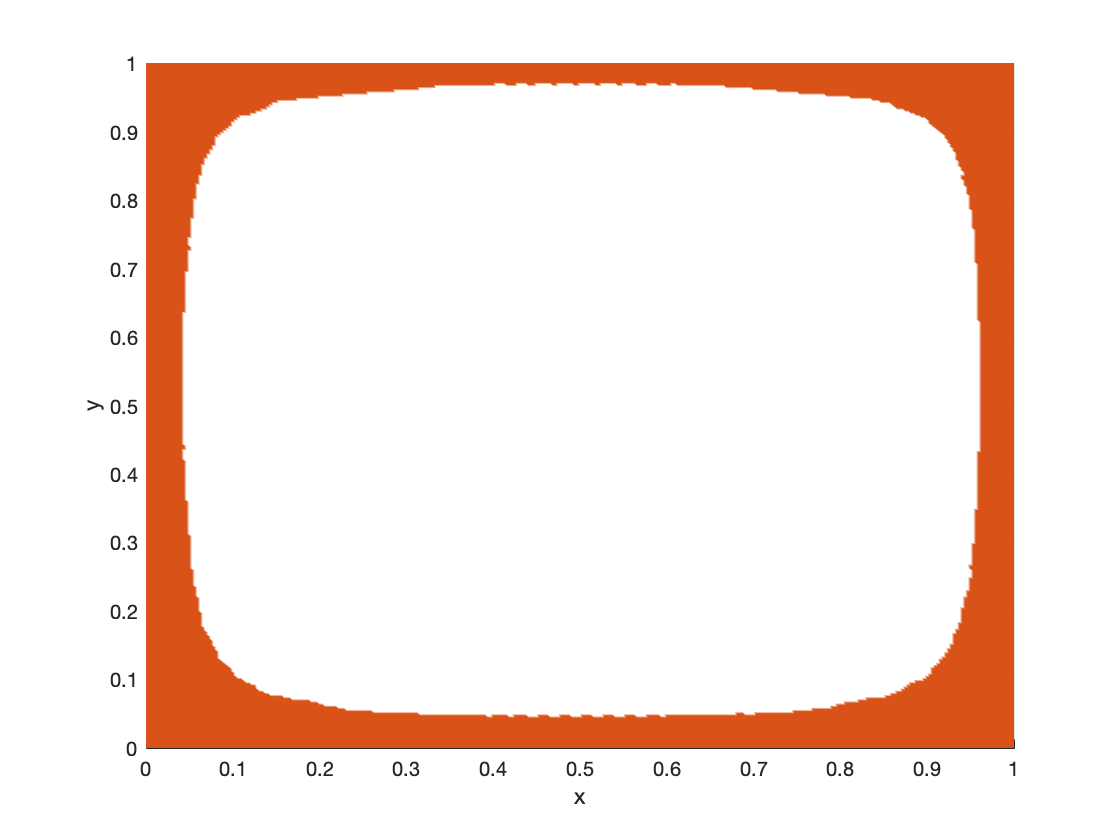}
\hspace{20pt}
 \includegraphics[width=0.3\linewidth]{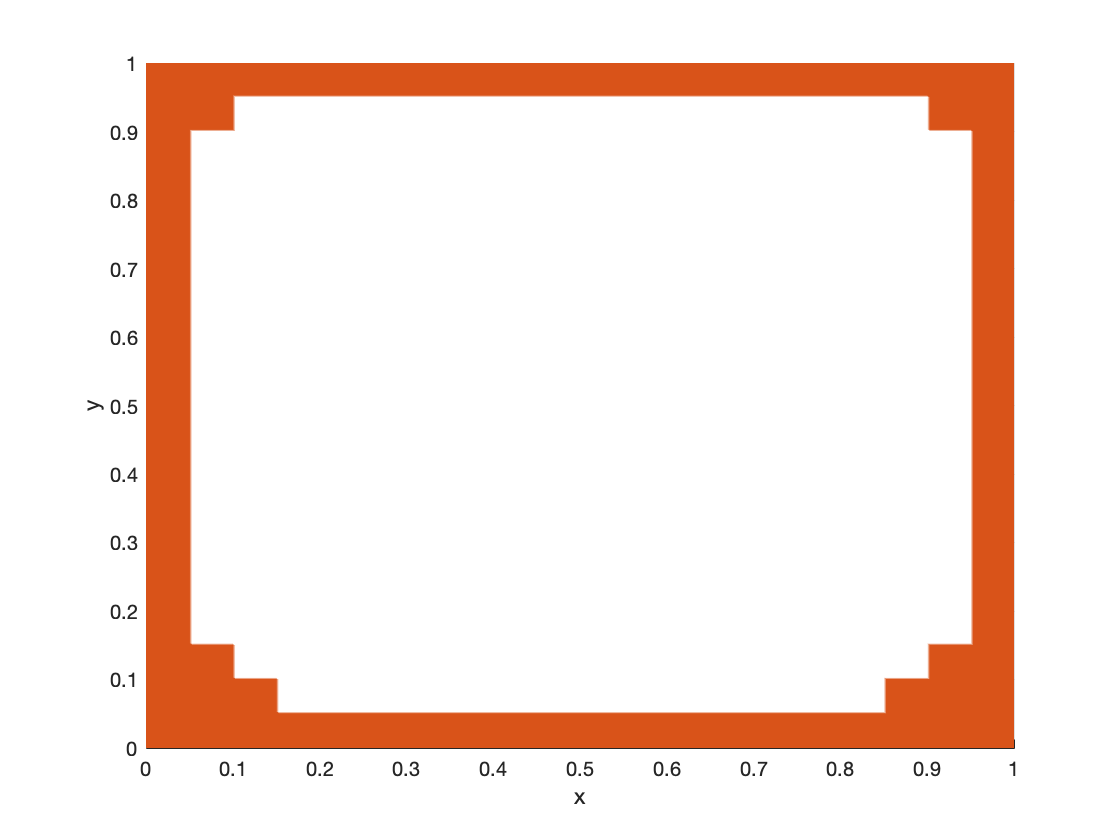}
\caption{Active set for $\phi_2$ for Example~\ref{example:HO}: $\bar u_h$ (left) and
$\buH$ (right, $H=1/20$).}
\label{fig:HOActiveSet2}
\end{figure}
\par\goodbreak
 The computation of the fine scale standard finite element solution
 $(\bar y_h,\bar u_h)$ of the discrete optimization problem takes
 $4.36\times 10^{+1}$
  seconds by using the PETSc/TAO library  with 20 processors.
 The computational time (in seconds) for $(\bMSLy,\buH)$
 using MATLAB on a laptop
 are presented in Table~\ref{table:HOTimes} for
 $H=1/10,1/20,1/40$.
 For $H=1/20$, the DD-LOD solution $(\MSLy,\buH)$ is a reasonable approximation of
 $(\bar y_h,\bar u_h)$
 (cf. Figures~\ref{fig:HOStateComparison}--\ref{fig:HOActiveSet2}) and its
 computation is more than 200 times faster than the computation of $(\bar y_h,\bar u_h)$.
\begin{table}[hhh!]
\centering
\begin{tabular}{|c|c|}
  \hline
  $H$ & Time\\
  \hline&\\[-12pt]
  $1/10$ & $1.71\times 10^{-2}$\\
  $1/20$ & $1.27\times 10^{-1}$\\
  $1/40$ & $1.40\times 10^{+1}$\\
  \hline
\end{tabular}
\par\medskip
\caption{Computational time in seconds for $(\bMSLy,\buH)$ for Example~\ref{example:HO}.}
\label{table:HOTimes}
\end{table}
\par
 In the second set of experiments we take $H=1/20$ and test the improved approximation by
 the DD-LOD solution $(\bMSLy,\bar u_\rho)$ for $\rho=1/40,1/80/1/160$ that is predicted
 by the estimates in Theorem~\ref{thm:DDLODError} and Theorem~\ref{thm:EnergyErrors}.
 This improvement can be observed by comparing the values of the cost function
 $\tilde J$ in Table~\ref{table:HOrhoMinimum} with the value
 $\tilde J(\bar y_h,\bar u_h)=-8.29631\times 10^{-5}$ for
 the fine scale solution.  The number of significant digits improves from 2 to $3$ as
 $\rho$ decreases from $1/20$ to $1/160$.
\par
\begin{table}[hhh!]
\centering
\begin{tabular}{|c | c|}
\hline &\\[-12pt]
$\rho$& $\tilde J(\bMSLy,\bar u_\rho)$\\
\hline &\\[-12pt]
$1/20$&$-8.28313\times10^{-5}$\\
$1/40$&$-8.29252\times10^{-5}$\\
$1/80$&$-8.29448\times10^{-5}$\\
$1/160$&$-8.29510\times 10^{-5}$\\
 \hline
\end{tabular}
\par\medskip
\caption{Values of $\tilde J(\bMSLy,\bar u_\rho)$ with $H=1/20$
 and various $\rho$ for Example~\ref{example:HO}.}
\label{table:HOrhoMinimum}
\end{table}
\par
 The improvement can also be visualized through a comparison of the graphs of the
 fine scale solution $\bar u_h$ and the DD-LOD solution $\bar u_\rho$ (with $H=1/20$ and $\rho=1/160$)
  in Figure~\ref{fig:HOrhoControlComparison}.  They are almost identical, which
 is not the case in Figure~~\ref{fig:HOControlComparison}.
\begin{figure}[hhh!]
\centering
  \includegraphics[width=0.3\linewidth]{OptControl-STD-319-osc-AL}
\hspace{20pt}
  \includegraphics[width=0.3\linewidth]{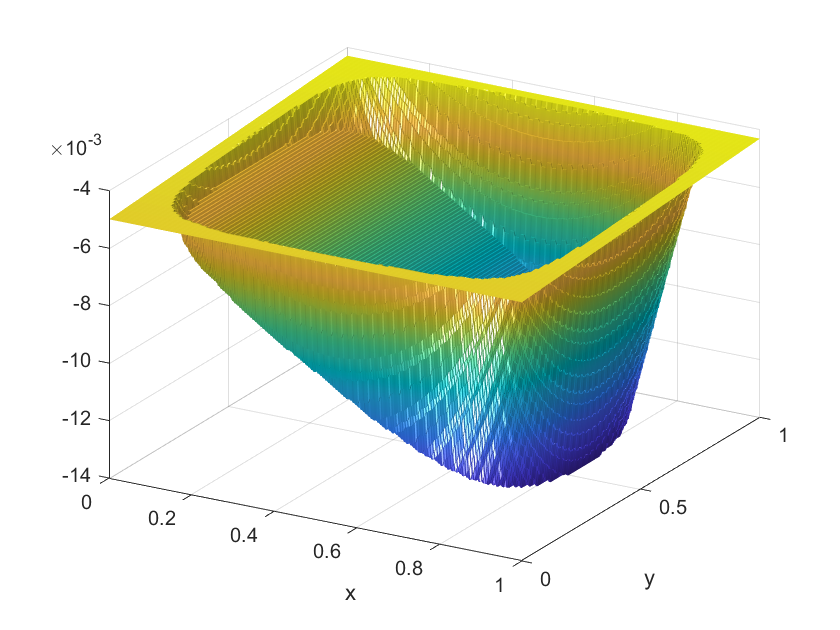}
\caption{Graph of $\bar u_h$ (left) and graph of $\bar u_\rho$ (right, $H=1/20$, $\rho=1/160$)
for Example~\ref{example:HO}.}
\label{fig:HOrhoControlComparison}
\end{figure}
\par
 We can also observe the improvement due to smaller  $\rho$ by comparing the active sets
 depicted in Figure~\ref{fig:HOrhoActiveSet1} and Figure~\ref{fig:HOrhoActiveSet2}.
 These sets are almost identical, which is not the case in
 Figure~\ref{fig:HOActiveSet1} and Figure~\ref{fig:HOActiveSet2}.
\begin{figure}[hhh!]
\centering
  \includegraphics[width=0.3\linewidth]{Active-set-part1-af-19-319-osc}
\hspace{20pt}
  \includegraphics[width=0.3\linewidth]{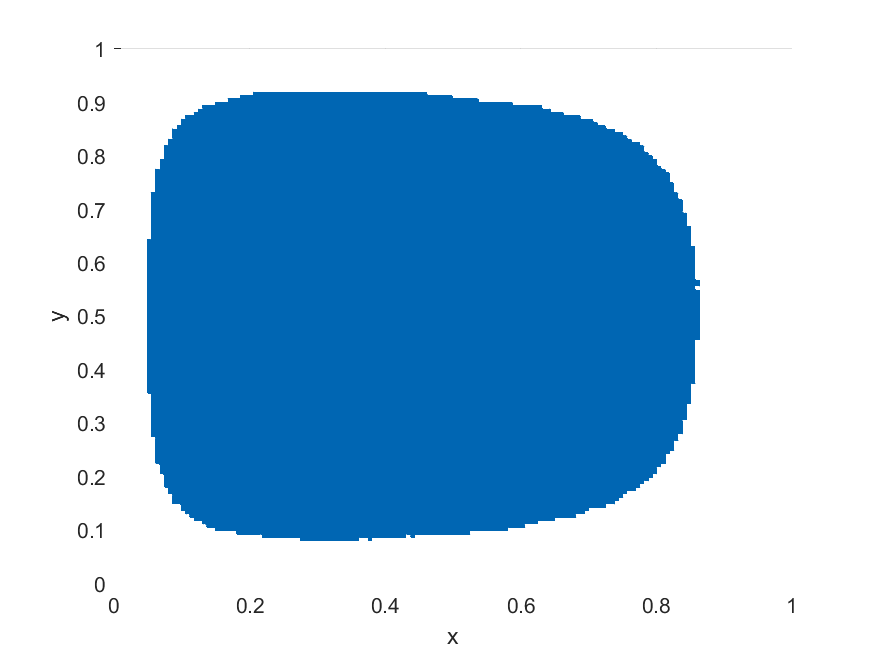}
\caption{Active set for $\phi_1$ for Example~\ref{example:HO}: $\bar u_h$ (left) and
$\bar u_\rho$ (right, $H=1/20$ and $\rho=1/160$).}
\label{fig:HOrhoActiveSet1}
\end{figure}
\begin{figure}[htb!]
\centering
  \includegraphics[width=0.3\linewidth]{Active-set-part2-af-19-319-osc}
\hspace{20pt}
 \includegraphics[width=0.3\linewidth]{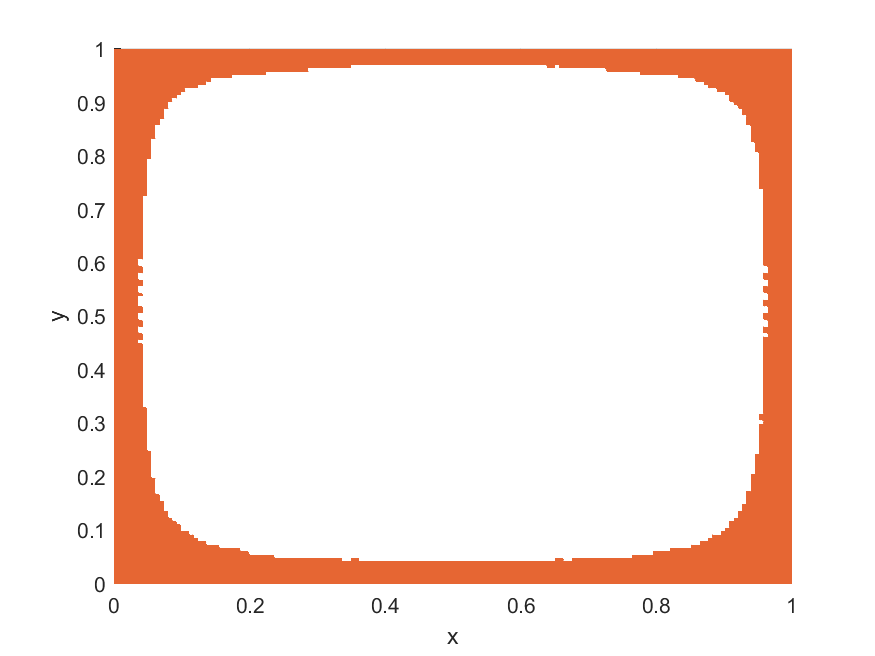}
\caption{Active set for $\phi_2$ for Example~\ref{example:HO}: $\bar u_h$ (left) and
$\bar u_\rho$ (right, $H=1/20$ and $\rho=1/160$).}
\label{fig:HOrhoActiveSet2}
\end{figure}
\end{example}
%%%%%%%%%%%%%%%%%%%%%%%
\section{Concluding Remarks}\label{sec:Conclusions}
 We have constructed and analyzed a multiscale finite element method for the
 optimal control problem defined by \eqref{eq:OCP}--\eqref{eq:aDef}.
 We showed that the approximate solution obtained by the DD-LOD finite element method
 on the coarse mesh $\cT_H$
 is, up to an $O(H^2+\rho)$ term for the $L_2$ error and an $O(H+\rho)$ term for
 the energy error,
  as good as the approximate solution obtained by
 a standard finite
 element method on a fine mesh $\cT_h$.   Alternatively we can say that up to the
 fine scale error
  the performance of the DD-LOD method is as good as standard finite element
  methods for smooth problems.
\par
 The DD-LOD multiscale finite element method is one of the simplest
 multiscale finite element methods in terms of construction and analysis.
 There is inherent parallelism in the construction of the DD-LOD finite
 element space that comes from domain
 decomposition so that it can readily benefit from high performance computing
  (cf. \cite{BGS:2022:LOD_OCP}), and its analysis
  only requires basic knowledge in finite
 element methods, domain decomposition
 methods and numerical linear algebra.
 After a multiscale basis has been
 computed off-line, the
 on-line solution with the coarse scale DD-LOD finite element method is fast.
 The multiscale finite element method is particularly useful for applications where
 the optimal control problem has to be solved repeatedly for different
 $y_d$, $\phi_1$ and $\phi_2$.
\par
 We note that the error estimates in Theorem~\ref{thm:AbstractErrorEstimate}
 and Theorem~\ref{thm:AbstractEnergyError} are applicable
 to any subspace $V_*$ of $H^1_0(\O)$.   The key is to have good error estimates for the
  Galerkin solution of
 \eqref{eq:RBVP}. In particular, we can take $V_*$ to be the
 LOD multiscale finite element spaces
 in \cite{PS:2016:Contrast,HM:2017:Contrast,HP:2013:LOD}
  and arrive at similar results.  Note that the LOD methods in
  \cite{PS:2016:Contrast,HM:2017:Contrast} are suitable for problems with high contrast.
\par
 We can also take $V_*$ to be the multiscale finite element space $V_h$
 from \cite{HW:1999:MS,HWC:1999:Multiscale,EH:2009:MSFEM} for problems
 with highly oscillatory and periodic coefficients
 (such as the problem in Example~\ref{example:HO}), where $h$ stands for the coarse
 mesh size.  The corresponding $L_2$ error estimate then takes the form
\begin{equation*}
  \|\bar y-\byh\|_\LT+\|\bar u-\buh\|_\LT+\|\bar p-\bph\|_\LT
  \leq C\Big(h^2+\epsilon+\frac{\epsilon}{h}+\rho\Big),
\end{equation*}
 where $\epsilon \,(<h)$ is the parameter for the small scale, and the positive
 constant $C$ only depends on
  $\|y_d\|_\LT$, $\|\phi_1\|_{H^1(\O)}$, $\|\phi_2\|_{H^1(\O)}$,
  $\gamma^{-1}$, $\alpha^{-1}$
 and the shape regularities of $\cT_h$ and $\cT_\rho$.
\par
 Similarly, the multiscale finite element methods in
 \cite{AC:2022:DD26,CLZZ:2023:MultiscaleCC}
 can also be analyzed by
 Theorem~\ref{thm:AbstractErrorEstimate}, Theorem~\ref{thm:AbstractEnergyError}
  and the estimates in
 \cite{OZB:2014:Homogenization,CEL:2018:Multiscale}.

%%%%%%%%%%%%%%%%%%%%%%%
%%%%%%%%%%%%%%%%%%%%%%%%%%%%%%%%%%%%%%
\section*{Acknowledgements}
Portions of this research were conducted with high performance computing
resources provided by Louisiana State University (http://www.hpc.lsu.edu).
%%%%%%%%%%%%%%%%%%%%%%%%%%%%%%%%%%%%%%
\section*{Funding}
This work was supported in part
 by the National Science Foundation under Grant No.
 DMS-19-13035 and Grant No. DMS-22-08404.
%%%%%%%%%%%%%%%%%%%%%%%%%%%%%%%%%%%%%%
\section*{Data Availability}
The datasets generated during and/or analyzed during the current study
are available from the corresponding author on reasonable request.
%%%%%%%%%%%%%%%%%%%%%%%

\end{document}